\documentclass[secthm,seceqn,amsthm]{elsart}
\usepackage{amsthm}		%for defining new theorem styles
\usepackage{amsmath}		%for subequations, align, multline environments
\usepackage{amsfonts}		%for Fraktur font and blackboard bold
\usepackage{amssymb}		%for 'rtimes' symbol, etc.

\theoremstyle{remark}
\newtheorem{remarkaftertheorem}{Remark}[thm]
\newtheorem{exampleaftertheorem}[remarkaftertheorem]{Example}


\journal{J.~Differential Equations}

\newcommand{\defeq}{\stackrel{\rm{def}}{=}}
\newcommand{\Hl}{\mathop{{}\it Hl}\nolimits}
\def\tfrac#1#2{{\textstyle\frac{#1}{#2}}}
\def\ri{{\rm i}}

\makeatletter
\newenvironment{sizeequation}[1]{%
  \skip@=\baselineskip
  #1%
  \baselineskip=\skip@
  \equation
}{\endequation \ignorespacesafterend}
\makeatother

\let\NoHyper\relax
\begin{document}
\begin{frontmatter}
\title{On Reducing the Heun Equation to the Hypergeometric Equation}
\author{Robert S. Maier\thanksref{NSF}}
\ead{rsm@math.arizona.edu}
\ead[url]{http://www.math.arizona.edu/\~{}rsm}
\address{Depts.\ of Mathematics and Physics, University of Arizona, Tucson AZ 85721, USA}
\thanks[NSF]{Partially supported by NSF grant PHY-0099484.}
\begin{abstract}
The reductions of the Heun equation to the hypergeometric equation by
polynomial transformations of its independent variable are enumerated and
classified.  Heun-to-hypergeometric reductions are similar to classical
hypergeometric identities, but the conditions for the existence of a
reduction involve features of the Heun equation that the hypergeometric
equation does not possess; namely, its cross-ratio and accessory
parameters.  The reductions include quadratic and cubic transformations,
which may be performed only~if the singular points of the Heun equation
form a harmonic or an equianharmonic quadruple, respectively; and several
higher-degree transformations.  This result cor\nobreak{rects} and extends
a theorem in a previous paper, which found only the quadratic
transformations.  [See K.~Kuiken, ``Heun's equation and the hypergeometric
equation'', {\em SIAM Journal on Mathematical Analysis\/} 10 (3) (1979),
655--657.]
\end{abstract}
\end{frontmatter}

\section{Introduction}
\label{sec:intro}
Consider the class of linear second-order differential equations on the
Riemann sphere~$\mathbb{CP}^1$ which are Fuchsian, i.e., have only regular
singular points~\cite{Hille76}.  Any such equation with exactly three
singular points can be transformed to the hypergeometric equation by
appropriate changes of the independent and dependent variables.  Similarly,
any such equation with exactly four singular points can be transformed to
the Heun equation.  (See
\cite[Chapter~15]{Erdelyi53},\cite{Ronveaux95,Snow52}.)

Solutions of the Heun equation are much less well understood than
hypergeometric functions~\cite{Arscott81}.  No~general integral
representation for them is known, for instance.  Such solutions have
recently been used in fluid dynamics~\cite{Craster98,Schmitz94} and
drift--diffusion theory~\cite{Debosscher98}.  They also arise in lattice
combinatorics~\cite{Guttmann93,Joyce94}.  But it~is difficult to carry~out
practical computations involving them.  An~explicit solution to the
two-point connection problem for the general Heun equation is not
known~\cite{Schafke80a}, though the corresponding problem for the
hypergeometric equation has a classical solution.  Most work on solutions
of the Heun equation has focused on special cases, such as the Lam\'e
equation~\cite{Erdelyi53,Maier04}.

Determining which Heun equation solutions are expressible in~terms of more
familiar functions would obviously be useful: it~would facilitate the
solution of the two-point connection problem, and the computation of Heun
equation monodromies.  A~significant result in this direction was obtained
by Kuiken~\cite{Kuiken79}.  It~is sometimes possible, by performing a
quadratic change of the independent variable, to reduce the Heun equation
to the hypergeometric equation, and thereby express its solutions in~terms
of hypergeometric functions.  Kuiken's quadratic transformations are not so
well known as they should~be.  The useful monograph edited by
Ronveaux~\cite{Ronveaux95} does not mention them explicitly, though it
lists Ref.~\cite{Kuiken79} in its bibliography.  One of Kuiken's
transformations was recently rediscovered by Ivanov~\cite{Ivanov2001}, in a
disguised form.

Unfortunately, the main theorem of Ref.~\cite{Kuiken79} is incorrect.  The
theorem asserts that a reduction to the hypergeometric equation, by a
rational change of the independent variable, is possible only~if the
singular points of the Heun equation form a harmonic quadruple in the sense
of projective geometry; in~which case the change of variables must be
quadratic.  In~this paper, we show that there are many alternatives.
A~reduction may also be possible if the singular points form an
equianharmonic quadruple, with the change of variables being cubic.
Additional singular point configurations permit changes of variable of
degrees $3$,~$4$, $5$, and~$6$.  Our main theorem (Theorem~\ref{thm:main})
and its corollaries classify all such reductions, up~to affine
automorphisms of the Heun and hypergeometric equations.  It replaces the
theorem of Ref.~\cite{Kuiken79}.

It follows from Theorem~\ref{thm:main} that in a suitably defined
`nontrivial' case, the local Heun function~$\Hl$ can be reduced to the
Gauss hypergeometric function~${}_2F_1$ by a formula of the type
$\Hl(t)={}_2F_1(R(t))$ only~if the pair $(d,q/\alpha\beta)$, computed from
the parameters of~$\Hl$, takes one of exactly $23$~values.  These are
listed in Theorem~\ref{thm:culmination}.  A~representative list of
reductions is given in Theorem~\ref{thm:useful0}.  These theorems should be
of~interest to special function theorists and applied mathematicians.
We~were led to our correction and expansion of the theorem of
Ref.~\cite{Kuiken79} by a discovery of Clarkson and
Olver~\cite{Clarkson96}: an~unexpected reduction of the Weierstrass form of
the equianharmonic Lam\'e equation to the hypergeometric equation.
In~\S\,\ref{sec:CO}, we explain how this is a special case of the cubic
Heun-to-hypergeometric reduction.

The new reductions are similar to classical hypergeometric transformations.
(See~\cite[Chapter~3]{Andrews99}; also~\cite[Chapter~2]{Erdelyi53}.)  But
reducing the Heun equation to the hypergeometric equation is more difficult
than transforming the hypergeometric equation to itself, since conditions
involving its singular point location parameter and accessory parameter,
as~well as its exponent parameters, must be satisfied.  Actually, the
reductions classified in this paper are of a somewhat restricted type,
since unlike many classical hypergeometric transformations, they involve
no~change of the dependent variable.  A~classification of reductions of the
more general type is possible, but is best phrased in algebraic-geometric
terms, as a classification of certain branched covers of the Riemann sphere
by itself.  A~further extension would allow the transformation of the
independent variable to be algebraic rather than polynomial or rational,
since at~least one algebraic Heun-to-hypergeometric reduction is known to
exist~\cite{Joyce94}.  Extended classification schemes are deferred to one
or more further papers.

\section{Preliminaries}
\subsection{The Equations}
\label{sec:defs}
The Gauss hypergeometric equation is
\def\theequation{$\mathfrak{h}$}
\begin{equation}
\label{eq:hyper}
\frac{\d^2 y}{\d z^2} + \left(\frac{c}{z} + \frac{a+b-c+1}{z-1}\right)
\frac{\d y}{\d z} + \frac{ab}{z(z-1)}\,y = 0,
\end{equation}\addtocounter{equation}{-1}\def\theequation{\arabic{section}.\arabic{equation}}where
$a,b,c\in\mathbb{C}$ are parameters.  It~and its solution space are
specified by the Riemann $P$-symbol
\begin{sizeequation}
{\small}
\label{eq:HeunP}
P\left\{
\begin{array}{cccc}
0&1&\infty& \\
0&0&a& ;z \\
1-c&c-a-b&b& 
\end{array}
\right\},
\end{sizeequation}
where each column, except the last, refers to a regular singular point.
The first entry is its location, and the final two are the characteristic
exponents of the solutions there.  The exponents at each singular point are
obtained by solving an indicial equation~\cite{Hille76}.  In~general, each
finite singular point~$z_0$ has $\zeta$~as an exponent if~and only~if the
equation has a local (Frobenius) solution of the form $(z-z_0)^\zeta h(z)$
in a neighborhood of~$z=z_0$, where $h$~is analytic and nonzero at~$z=z_0$.
If~the exponents at~$z=z_0$ differ by an integer, this statement must be
modified: the solution corresponding to the smaller exponent may have a
logarithmic singularity at~$z=z_0$.  The definition extends in a
straightforward way to~$z_0=\infty$, and also to ordinary points, each of
which has exponents~$0,1$.

There are $2\times3=6$ local solutions of~($\mathfrak{h}$) in~all: two per
singular point.  If $c$~is not a nonpositive integer, the solution at~$z=0$
belonging to the exponent zero will be analytic.  When normalized to unity
at~$z=0$, it~will be the Gauss hypergeometric function
${}_2F_1(a,b;c;z)$~\cite{Erdelyi53}.  This is the sum of a hypergeometric
series, which converges in a neighborhood of~$z=0$.  In~general,
${}_2F_1(a,b;c;z)$ is not defined when $c$~is a nonpositive integer.

The Heun equation is usually written in the form
\def\theequation{$\mathfrak{H}$}
\begin{equation}
\label{eq:Heun}
\frac{\d^2 u}{\d t^2}
+ \left( \frac\gamma t + \frac\delta{t-1} + \frac\epsilon{t-d}
  \right)\frac{\d u}{\d t} + \frac{\alpha\beta t - q}{t(t-1)(t-d)}\,u = 0.
\end{equation}\addtocounter{equation}{-1}\def\theequation{\arabic{section}.\arabic{equation}}Here $d\in\mathbb{C}$, the location of the fourth singular point,
is a parameter~($d\neq0,1$), and
$\alpha,\beta,\gamma,\delta,\epsilon\in\mathbb{C}$ are exponent-related
parameters.  The $P$-symbol is
\begin{sizeequation}
{\small}
\label{eq:Psymbol}
P\left\{
\begin{array}{ccccc}
0&1&d&\infty& \\
0&0&0&\alpha& ;t \\
1-\gamma&1-\delta&1-\epsilon&\beta& 
\end{array}
\right\}.
\end{sizeequation}
This does not uniquely specify the equation and its solutions, since it
omits the accessory parameter~$q\in\mathbb{C}$.  The exponents are
constrained by
\begin{equation}
\label{eq:Pconstraint}
\alpha+\beta-\gamma-\delta-\epsilon+1 = 0.
\end{equation}
This is a special case of Fuchs's relation, according to which the sum of
the $2n$~characteristic exponents of any second-order Fuchsian equation
on~$\mathbb{CP}^1$ with $n$~singular points must equal
$n-2$~\cite{Poole36}.

There are $2\times4=8$ local solutions of~($\mathfrak{H}$) in~all: two per
singular point.  If $\gamma$~is not a nonpositive integer, the solution
at~$t=0$ belonging to the exponent zero will be analytic.  When normalized
to unity at~$t=0$, it~is called the local Heun function, and is denoted
$\Hl(d,q;\alpha,\beta,\gamma,\delta;t)$~\cite{Ronveaux95}.  It~is the sum
of a Heun series, which converges in a neighborhood
of~$t=0$~\cite{Ronveaux95,Snow52}.  In~general,
$\Hl(d,q;\alpha,\beta,\gamma,\delta;t)$ is not defined when $\gamma$~is a
nonpositive integer.

If $\epsilon=0$ and $q=\alpha\beta d$, the Heun equation loses a singular
point and becomes a hypergeometric equation.  Similar losses occur if
$\delta=0$, $q=\alpha\beta$, or $\gamma=0$,~$q=0$.  This paper will exclude
the case when the Heun equation has fewer than four singular points, since
reducing ($\mathfrak{h}$) to itself is a separate problem, leading to the
classical hypergeometric transformations.  The following case, in which the
solution of~(\ref{eq:Heun}) can be reduced to quadratures, will be
initially excluded.
\begin{defn}
\label{def:triviality}
If $\alpha\beta=0$ and~$q=0$, the Heun equation~(\ref{eq:Heun}) is said to
be trivial.  Triviality implies that one of the exponents at~$t=\infty$ is
zero (i.e., $\alpha\beta=0$), and is implied by absence of the singular
point at~$t=\infty$ (i.e., $\alpha\beta=0$, $\alpha+\beta=1$, $q=0$).
\end{defn}

The transformation to ($\mathfrak{H}$) or~($\mathfrak{h}$) of a linear
second-order Fuchsian differential equation with singular points at
$t=0,1,d,\infty$ (resp.\ $z=0,1,\infty$), and with arbitrary exponents, is
accomplished by certain linear changes of the dependent variable, called
F-homotopies. (See~\cite{Erdelyi53} and~\cite[\S\,{A}2 and
Addendum,~\S\,1.8]{Ronveaux95}.)  If~an equation with singular points
at~$t=0,1,d,\infty$ has dependent variable~$u$, carrying~out the
substitution $\tilde u(t)=t^{-\rho}(t-1)^{-\sigma}(t-d)^{-\tau} u(t)$ will
convert the equation to a new one, with the exponents at~$t=0,1,d$ reduced
by~$\rho,\sigma,\tau$ respectively, and those at~$t=\infty$ increased by
$\rho+\sigma+\tau$.  By~this technique, one exponent at each finite
singular point can be shifted to zero.

In fact, the Heun equation has a group of F-homotopic automorphisms
isomorphic to~$({\mathbb Z}_2)^3$, since at each of $t=0,1,d$, the
exponents~$0,\zeta$ can be shifted to~$-\zeta,0$, i.e., to~$0,-\zeta$.
Similarly, the hypergeometric equation has a group of F-homotopic
automorphisms isomorphic to $({\mathbb Z}_2)^2$.  These groups act on the
$6$~and~$3$-dimensional parameter spaces, respectively.  For example, one
of the latter actions is $(a,b;c)\mapsto(c-a,c-b;c)$, which is induced by
an F-homotopy at~$z=1$.  From this F-homotopy follows Euler's
transformation~\cite[\S\,2.2]{Andrews99}
\begin{equation}
\label{eq:flip}
{}_2F_1(a,\,b;\,c;\,z)= (1-z)^{c-a-b}{}_2F_1(c-a,\,c-b;\,c;\,z),
\end{equation}
which holds because ${}_2F_1$~is a local solution at~$z=0$, rather than
at~$z=1$.

If the singular points of the differential equation are arbitrarily placed,
transforming it to the Heun or hypergeometric equation will require a
M\"obius (i.e., projective linear or homographic) transformation, which
repositions the singular points to the standard locations.  A~unique
M\"obius transformation maps any three distinct points in~$\mathbb{CP}^1$
to any other three; but the same is not true of four points, which is why
($\mathfrak{H}$)~has the singular point~$d$ as a free parameter.

\subsection{The Cross-Ratio}
\label{subsec:crossratio}
The characterization of Heun equations that can be reduced to the
hypergeometric equation will employ the cross-ratio orbit of
$\{0,1,d,\infty\}$, defined as follows.  If $A,B,C,D\in\mathbb{CP}^1$ are
distinct, their cross-ratio is
\begin{equation}
(A,B;C,D)\defeq
\frac{(C-A)(D-B)}{(D-A)(C-B)}\in\mathbb{CP}^1\setminus\{0,1,\infty\},
\end{equation}
which is invariant under M\"obius transformations.  Permuting $A,B,C,D$
yields an action of the symmetric group~$S_4$
on~$\mathbb{CP}^1\setminus\{0,1,\infty\}$.  The cross-ratio is invariant
under interchange of $A,B$ and~$C,D$, and also under simultaneous
interchange of the two points in each pair.  So each orbit contains no~more
than $4!/4=6$ cross-ratios.  The possible actions of~$S_4$
on~$s\in\mathbb{CP}^1\setminus\{0,1,\infty\}$ are generated by
$s\mapsto1-s$ and $s\mapsto 1/s$, and the orbit of~$s$ comprises
\begin{equation}
s,\quad 1-s,\quad 1/s,\quad 1/(1-s),\quad s/(s-1),\quad (s-1)/s,
\end{equation}
which may not be distinct.  This is called the cross-ratio orbit of~$s$;
or, if~$s=(A,B;\allowbreak C,D)$, the cross-ratio orbit of the unordered
set $\{A,B,C,D\}\subset\mathbb{CP}^1$.  Two sets of distinct points
$\{A_i,B_i,C_i,D_i\}$ ($i=1,2$) have the same cross-ratio orbit iff they
are related by a M\"obius transformation.

Cross-ratio orbits generically contain six values, but there are two
exceptional orbits: one with three and one with two.  If $(A,B;\allowbreak
C,D)=-1$, the cross-ratio orbit of $\{A,B,C,D\}$ will be
$\{-1,\frac12,2\}$.  The value~$-1$ for $(A,B;\allowbreak C,D)$ defines a
so-called harmonic configuration: $A,B$ and $C,D$ are said to be harmonic
pairs.  More generally, if $\{A,B,C,D\}$ has cross-ratio orbit
$\{-1,\frac12,2\}$, it~is said to be a harmonic quadruple.  It~is easy to
see that if $C=\infty$ and $A,B,D$ are distinct finite points, then $A,B$
and $C,D$ will be harmonic pairs iff $D$~is the midpoint of the line
segment~$\overline{AB}$.  In~consequence, $\{A,B,\infty,D\}$ will be a
harmonic quadruple iff $\{A,B,D\}\subset\mathbb{C}$ comprises collinear,
equally spaced points.  So, $\{A,B,C,D\}\subset\mathbb{CP}^1$ will be a
harmonic quadruple iff it can be mapped by a M\"obius transformation to a
set consisting of three equally spaced finite points and the point at
infinity; equivalently, to the vertices of a square in~$\mathbb{C}$.

The cross-ratio orbit containing exactly two values is
$\{\frac12\pm\ri\frac{\sqrt3}2\}$.  Any set $\{A,B,C,D\}$ with this as
cross-ratio orbit is said to be an equianharmonic quadruple.
$\{A,B,\infty,D\}$ will be an equianharmonic quadruple iff $A,B,D$ are the
vertices of an equilateral triangle in~$\mathbb{C}$\null.  If
$\mathbb{CP}^1$ is interpreted as a sphere via the usual stereographic
projection, then by an affine transformation (a~special M\"obius
transformation), this situation reduces to the case when $A,B,\infty,D$ are
the vertices of a regular tetrahedron.  So,
$\{A,B,C,D\}\subset\mathbb{CP}^1$ will be an equianharmonic quadruple iff
it can be mapped by a M\"obius transformation to the vertices of a regular
tetrahedron in~$\mathbb{CP}^1$.

Cross-ratio orbits are of two sorts: real orbits such as the harmonic
orbit, and non-real orbits such as the equianharmonic orbit.  All values in
a real orbit are real, and in a non-real orbit, all have a nonzero
imaginary part.  So, $\{A,B,C,D\}$ will have a specified real orbit as its
cross-ratio orbit iff it can be mapped by a M\"obius transformation to a
set consisting of three specified collinear points in~$\mathbb{C}$ and the
point at infinity; equivalently, to the vertices of a specified quadrangle
(generically irregular) in~$\mathbb{C}$\null.  Similarly, it will have a
specified non-real orbit as its cross-ratio orbit iff it can be mapped to a
set consisting of a specified triangle in~$\mathbb{C}$ and the point at
infinity; equivalently, to the vertices of a specified tetrahedron
(generically irregular) in~$\mathbb{CP}^1$.

The cross-ratio orbit of $\{0,1,d,\infty\}$ will be the harmonic orbit iff
$d$~equals $-1$, $\frac12$, or~$2$, and the equianharmonic orbit iff $d$~equals
$\frac12\pm \ri\frac{\sqrt{3}}2$.  In~contrast, it will be a specified
generic orbit iff $d$~takes one of six orbit-specific values.
The cross-ratio orbit of $\{0,1,d,\infty\}$ being a specified orbit is
equivalent to its being the same as the cross-ratio orbit of some specified
quadruple of the form $\{0,1,D,\infty\}$, i.e., to there being a M\"obius
transformation that maps $\{0,1,d,\infty\}$ onto $\{0,1,D,\infty\}$.
The possibilities are
\begin{equation}
\label{eq:firstlist}
d=D,\quad 1-D,\quad 1/D,\quad 1/(1-D),\quad D/(D-1),\quad (D-1)/D.
\end{equation}
By~examination, this is equivalent to $\{0,1,d\}$ being mapped onto
$\{0,1,D\}$ by some {\em affine\/} transformation, i.e., to $\triangle01d$
being similar to~$\triangle01D$.  The corresponding affine transformations
$t\mapsto A_1(t)$ are
\begin{equation}
\label{eq:secondlist}
A_1(t)= t,\quad 1-t,\quad Dt,\quad (D-1)t+1,\quad (1-D)t+D,\quad D(1-t).
\end{equation}
This interpretation gives a geometric significance to the possible values
of~$d$.

For the equianharmonic orbit, in which the six values degenerate to two,
the triangle $\triangle01D$ may be taken to be any equilateral triangle.
For any real orbit, $\triangle01D$ must be taken to be degenerate, with
collinear vertices.  For example, the cross-ratio orbit of
$\{0,1,d,\infty\}$ being the harmonic orbit is equivalent to $\triangle01d$
being similar to~$\triangle012$, i.e., to $\{0,1,d\}$ consisting of three
equally spaced collinear points.  This will be the case iff $d$~equals
$-1$, $\frac12$, or~$2$, in agreement with the definition of a harmonic
quadruple.

\subsection{Automorphisms}
\label{subsec:auto}
According to the theory of the Riemann $P$-function, any M\"obius
transformation~$M$ of the independent variable will preserve characteristic
exponents.  For the hypergeometric equation~(\ref{eq:hyper}), this implies
that if $M$~is one of the $3!$~M\"obius transformations that permute the
singular points $z=0,1,\infty$, the exponents of the transformed
equation~($\tilde{\mathfrak{h}}$) at its singular points
$M(0),M(1),M(\infty)$ will be those of~(\ref{eq:hyper}) at~$0,1,\infty$.
But if $M$~is not affine, i.e., $M(\infty)\neq\infty$, then
($\tilde{\mathfrak{h}}$)~will not in general be a hypergeometric equation,
since its exponents at~$M(\infty)$ may both be nonzero.  To~convert
($\tilde{\mathfrak{h}}$) to a hypergeometric equation, the permutation must
in this case be followed by an F-homotopic transformation of the form
$\tilde y(z)=[z-M(\infty)]^{-a} y(z)$ or $\tilde y(z)=[z-M(\infty)]^{-b}
y(z)$.

\begin{defn}
${\rm Aut}(\mathfrak{h})$, the automorphism group of the hypergeometric
equation~(\ref{eq:hyper}), is the group of changes of variable (M\"obius of
the independent variable, linear of the dependent) which
leave~(\ref{eq:hyper}) invariant, up~to parameter changes.  Similarly,
${\rm Aut}(\mathfrak{H})$ is the automorphism group of~($\mathfrak{H}$).
\end{defn}

${\rm Aut}(\mathfrak{h})$ acts on the 3-dimensional parameter space
of~(\ref{eq:hyper}).  It~contains the symmetric group~$S_3$ of permutations
of singular points as a subgroup, and the group $({\mathbb Z}_2)^2$ of
F-homotopies as a normal subgroup.  So ${\rm Aut}(\mathfrak{h})\simeq
({\mathbb Z}_2)^2 \rtimes S_3$, a~semidirect product.  It~is isomorphic
to~$S_4$, the octahedral group~\cite{Dwork84}.

\begin{defn}
Within ${\rm Aut}(\mathfrak{h})$, the M\"obius automorphism subgroup is the
group ${\mathcal M}(\mathfrak{h})\defeq\{1\}\times S_3$, which permutes the
singular points $z=0,1,\infty$.  The subgroup of affine automorphisms is
${\mathcal A}(\mathfrak{h})\defeq\{1\}\times S_2$, which permutes the\/
{\em finite\/} singular points $z=0,1$, and fixes~$\infty$.  (It~is
generated by the involution $z\mapsto1-z$.)  The F-homotopic automorphism
subgroup is $({\mathbb Z}_2)^2\times\{1\}$.
\end{defn}

The action of ${\rm Aut}(\mathfrak{h})$ on the $2\times3=6$ local solutions
is as~follows. $\left|{\rm Aut}(\mathfrak{h})\right|=2^2\times3!=24$, and
applying the transformations in~${\rm Aut}(\mathfrak{h})$ to any single
local solution yields $24$~solutions of~(\ref{eq:hyper}).  Applying them
to~${}_2F_1$, for instance, yields the $24$~series solutions of
Kummer~\cite{Dwork84}.  However, the $24$~solutions split into six sets of
four, since for each singular point $z_0\in\{0,1,\infty\}$ there is a
subgroup of~${\rm Aut}(\mathfrak{h})$ of order~$2^1\times2!=4$, each
element of~which fixes~$z=z_0$ and performs no~F-homotopy there; so it
leaves each local solution at~$z=z_0$ invariant.

For example, the four transformations in the subgroup associated to~$z=0$
yield four equivalent expressions for ${}_2F_1(a,b;c;z)$; one of which is
${}_2F_1(a,b;c;z)$ itself, and another of which appears above
in~(\ref{eq:flip}).  The remaining two are
$(1-z)^{-a}{}_2F_1\left(a,c-b;c;z/(z-1)\right)$ and
$(1-z)^{-b}{}_2F_1\left(b,c-a;c;z/(z-1)\right)$.  The five remaining sets
of four are expressions for the five remaining local solutions.  One that
will play a role is the `second' local solution at~$z=0$, which belongs to
the exponent $1-c$.  One of the four expressions for~it, in~terms
of~${}_2F_1$, is~\cite{Erdelyi53}
\begin{equation}
\label{eq:tilde1}
\widetilde{{}_2F_1}(a,\,b;\,c;\,z)\defeq z^{1-c}{}_2F_1(a-c+1,\,b-c+1;\,2-c;\,z).
\end{equation}
$\widetilde{{}_2F_1}(a,b;c;z)$~is defined if $c\neq2,3,4,\dotsc$\ The
second local solution must be specified differently if $c=2,3,4,\dotsc$,
and also if~$c=1$, since in that case, $\widetilde{{}_2F_1}$ reduces
to~${{}_2F_1}$.  (Cf.~\cite[\S\,15.5]{Abramowitz65}.)  When
$\widetilde{{}_2F_1}$ is defined, it may be given a unique meaning by
choosing the principal branch of~$z^{1-c}$.

The automorphism group of the Heun equation is slightly more complicated.
There are $4!$~M\"obius transformations~$M$ that map the singular points
$t=0,1,d,\infty$ onto $t=0,1,d',\infty$, for
some~$d'\in\mathbb{CP}^1\setminus\{0,1,\infty\}$.  The possible~$d'$
constitute the cross-ratio orbit of~$\{0,1,d,\infty\}$.  Of~these
$4!$~transformations, $3!$~fix $t=\infty$, i.e., are affine.  All
values~$d'$ are obtained by affine transformations, i.e., a mapping is
possible iff $\triangle01d$ is similar to~$\triangle01d'$.  (Cf.~the
discussion in~\S\,\ref{subsec:crossratio}.)  If~$M$ is not affine, it must
be followed by an F-homotopic transformation of the form $\tilde
u(t)=[t-M(\infty)]^{-\alpha} u(t)$ or $\tilde u(t)=[t-M(\infty)]^{-\beta}
u(t)$.

${\rm Aut}(\mathfrak{H})$ acts on the 6-dimensional parameter space
of~(\ref{eq:Heun}).  It~contains the group~$S_4$ of singular point
permutations as a subgroup, and the group $({\mathbb Z}_2)^3$ of
F-homotopies as a normal subgroup.  So ${\rm Aut}(\mathfrak{H})\simeq
({\mathbb Z}_2)^3 \rtimes S_4$.  It turns~out to be isomorphic to the
Coxeter group~$\mathcal{D}_4$~\cite{Grove85}.

\begin{defn}
Within ${\rm Aut}(\mathfrak{H})$, the M\"obius automorphism subgroup is the
group ${\mathcal M}(\mathfrak{H})\defeq\{1\}\times S_4$, which maps between
sets of singular points of the form $\{0,1,d',\infty\}$.  The subgroup of
affine automorphisms is ${\mathcal A}(\mathfrak{H})\defeq\{1\}\times S_3$,
which maps between sets of\/ {\em finite\/} singular points of the form
$\{0,1,d'\}$, and fixes~$\infty$.  The F-homotopic automorphism subgroup is
$({\mathbb Z}_2)^3\times\{1\}$.
\end{defn}

The action of ${\rm Aut}(\mathfrak{H})$ on the $2\times4=8$ local solutions
is as~follows. $\left|{\rm Aut}(\mathfrak{H})\right|=2^3\times4!=192$, and
applying the transformations in~${\rm Aut}(\mathfrak{H})$ to any single
local solution yields $192$~solutions of~(\ref{eq:Heun}).  However, the
$192$~solutions split into eight sets of~$24$, since for each singular
point $t_0\in\{0,1,d,\infty\}$ there is a subgroup of~${\rm
Aut}(\mathfrak{H})$ of order~$2^2\times3!=24$, each element of~which
fixes~$t=t_0$ and performs no~F-homotopy there; so it leaves each local
solution at~$t=t_0$ invariant.  This statement must be interpreted with
care: selecting $t_0=d$ selects not a single singular point, but rather a
cross-ratio orbit.

The $24$~transformations in the subgroup associated to~$t_0=0$ yield
$23$~equivalent expressions for $\Hl(d,q;\alpha,\beta,\gamma,\delta;t)$,
one of which, the~only one with no F-homotopic prefactor, appears on the
right in the identity~\cite{Ronveaux95,Snow52}
\begin{equation}
\label{eq:newguy1}
\quad\Hl\left(d,\,q;\,\alpha,\,\beta,\,\gamma,\,\delta;\,t\right)
=\Hl\left(1/d,\,q/d;\,\alpha,\,\beta,\,\gamma,\,\alpha+\beta-\gamma-\delta+1;\,t/d\right).
\end{equation}
(The two sides are defined if $\gamma$ is not a nonpositive integer.)  The
remaining seven sets of~24 are expressions for the remaining seven local
solutions.  One that will play a role is the `second' solution at~$t=0$,
which belongs to the exponent $1-\gamma$.  One of the $24$~expressions
for~it, in~terms of~$\Hl$, is~\cite{Snow52}
\begin{equation}
\label{eq:newguy2}
\quad\widetilde\Hl(d,\,q;\,\alpha,\,\beta,\,\gamma,\,\delta;\,t)\defeq
t^{1-\gamma}
\Hl(d,\,\tilde q;\,\alpha-\gamma+1,\,\beta-\gamma+1,\,2-\gamma,\,\delta;\,t),
\end{equation}
where the transformed accessory parameter~$\tilde q$ equals
$q+(1-\gamma)(\epsilon+d\delta)$.  The quantity
$\widetilde{\Hl}(d,q;\alpha,\beta,\gamma,\delta;t)$ is defined if
$\gamma\neq2,3,4,\dotsc$\ The second local solution must be specified
differently if $\gamma=2,3,4,\dotsc$, and also if~$\gamma=1$, since in that
case, $\widetilde{\Hl}$ reduces to~$\Hl$.  When $\widetilde\Hl$ is defined,
it may be given a unique meaning by choosing the principal branch
of~$t^{1-\gamma}$.

In general, transformations in~${\rm Aut}(\mathfrak{H})$ will alter not
merely $d$~and the exponent parameters, but also the accessory
parameter~$q$.  This is illustrated by (\ref{eq:newguy1})
and~(\ref{eq:newguy2}).  The general transformation law of~$q$ is rather
complicated.  Partly for this reason, no~satisfactory list of the
192~solutions has appeared in~print.  The original paper of
Heun~\cite{Heun1889} tabulates 48 of the~192, but omits the value of~$q$ in
each.  His table also unfortunately contains numerous misprints and cannot
be used in practical applications~\cite[\S\,6.3]{Schmitz94}.  Incidentally,
one sometimes encounters a statement that there are only 96 distinct
solutions~\cite{Babister67,Exton93,Ronveaux95}.  This is true only if one
uses~(\ref{eq:newguy1}) to identify the 192 solutions in pairs.

\section{Polynomial Heun-to-Hypergeometric Reductions}
\label{sec:main}
We now state and prove Theorem~\ref{thm:main}, our corrected and expanded
version of the theorem of Kuiken~\cite{Kuiken79}.

The theorem will characterize when a homomorphism of rational substitution
type from the Heun equation~(\ref{eq:Heun}) to the hypergeometric
equation~(\ref{eq:hyper}) exists.  It~will list the possible substitutions,
up~to affine automorphisms of the two equations.  It~is really a
characterization of the ${\mathcal A}(\mathfrak{H})$-orbits that can be
mapped by homomorphisms of this type to ${\mathcal
A}(\mathfrak{h})$-orbits.  The possible substitutions, it turns~out, are
all polynomial.

For ease of understanding, the characterization will be concrete: it~will
state that $\triangle01d$ must be similar to one of five specified
triangles of the form $\triangle01D$.  By the remarks
in~\S\,\ref{subsec:crossratio}, similarity occurs iff $d$~belongs to the
cross-ratio orbit of~$D$, i.e., iff $D$~can be generated from~$d$ by
repeated application of $d\mapsto 1-d$ and $d\mapsto1/d$.  The two
exceptional cross-ratio orbits, namely $\{-1,\frac12,2\}$ (harmonic) and
$\{\frac12\pm\ri\frac{\sqrt3}2\}$ (equianharmonic), will play a prominent
role.  It~is worth noting that if ${\rm Re}\, D=\frac12$, the orbit of~$D$
is closed under complex conjugation.

For each value of~$D$, the polynomial map from $\mathbb{CP}^1\ni t$ to
$\mathbb{CP}^1\ni z$, which is denoted~$R$, will be given explicitly
when~$d=D$.  If $d$~is any other value on the cross-ratio orbit of~$D$, as
listed in~(\ref{eq:firstlist}), the polynomial map would be computed by
composing with the corresponding affine transformation $A_1$
of~$\mathbb{C}$ that takes $\triangle01d$ to~$\triangle01D$; which is
listed in~(\ref{eq:secondlist}).  So if $d\neq D$, statements in the
theorem dealing with singular points, characteristic exponents, and the
accessory parameter must be altered.  For example, case~1 of the theorem
refers to a distinguished singular point~$d_0$, the mandatory value of
which is given when~$d=D$.  If~$d\neq D$, its mandatory value would be
computed as the preimage of that point under~$A_1$.  Similarly, a statement
``the exponents of $t=0$ must be~$0,1/2$'', valid when~$d=D$, would be
interpreted if $d\neq D$ as ``the exponents of the preimage of $t=0$
under~$A_1$ must be~$0,1/2$''.  And a statement that $q/\alpha\beta$ must
take some value would be interpreted if $d\neq D$ as a statement that
$q/\alpha\beta$ must equal the preimage of that value under~$A_1$.

In the statement of the theorem and what follows, $S\defeq1-R$.

\begin{thm}
\label{thm:main}
A Heun equation~{\rm(}\ref{eq:Heun}{\rm)}, which has four singular points
and is nontrivial {\rm(}i.e., $\alpha\beta\neq0$ or~$q\neq0${\rm)}, can be
transformed to the hypergeometric equation~{\rm(}$\mathfrak{h}${\rm)} by a
rational substitution $z=R(t)$ if~and only~if $R$~is a polynomial,
$\alpha\beta\neq0$, and one of the following two conditions is satisfied.

% change the numbering of the top-level enumerate environment
% to an upright version of (1),(2),...
\def\labelenumi{{\rm (\theenumi)}}        \def\theenumi{\arabic{enumi}}

\begin{enumerate}
\item $\triangle01d$ is similar to $\triangle01D$, for one of the values
of~$D$ listed in subcases {\rm 1a}--{\rm 1c}; each of which is real, so the
triangle must be degenerate.  Also, the normalized accessory parameter
$q/\alpha\beta$ must equal one of\/ $0,1,d$, which may be denoted~$d_0$.
Each subcase lists the value of~$d_0$ when~$d=D$.
\item $\triangle01d$ is similar to $\triangle01D$, for one of the values
of~$D$ listed in subcases {\rm 2a}--{\rm 2d}; each of which is non-real and
has real part equal to~$\frac12$, so the triangle must be isosceles.  Each
subcase lists the value of~$q/\alpha\beta$ when $d=D$.
\end{enumerate}

Besides specifying $D$ and the value of~$q/\alpha\beta$ when $d=D$, each
subcase imposes restrictions on the characteristic exponent parameters at
the singular points\/ $0,1,d$. The subcases of the `real' case~{\rm 1} are the
following.

% change the numbering of the top-level enumerate environment
% to yield (1a),(1b),... instead of (1),(2),...
\def\labelenumi{{\em (\theenumi)}} \def\theenumi{1\alph{enumi}}

\begin{enumerate}
\item {\rm [}Harmonic {\rm(}equally spaced collinear points\,{\rm)}
case.{\rm]} $D=2$.  Suppose $d=D$.  Then $d_0$~must equal\/~$1$, and
$t=0,d$ must have the same characteristic exponents, i.e.,
$\gamma=\epsilon$.  In~general, either $R$ or~$S$ will be the degree-$2$
polynomial $t(2-t)$, which maps $t=0,d$ to~$z=0$ and $t=1$ to~$z=1$
{\rm(}with double multiplicity\/{\rm)}.  There are special circumstances in
which $R$~may be quartic, which are listed separately, as subcase~{\rm 1c}.

\item $D=4$.  Suppose $d=D$.  Then $d_0$~must equal\/~$1$, the point $t=1$
must have characteristic exponents that are double those of~$t=d$, i.e.,
$1-\delta=2(1-\epsilon)$, and $t=0$ must have exponents $0,1/2$, i.e.,
$\gamma=1/2$.  Either $R$ or~$S$ will be the degree-$3$ polynomial
$(t-1)^2(1-t/4)$, which maps $t=0$ to~$z=1$ and $t=1,d$ to~$z=0$ {\rm(}the
former with double multiplicity\/{\rm)}.

\item {\rm [}Special harmonic case.{\rm]} $D=2$.  Suppose $d=D$.  Then
$d_0$~must equal~$1$, and $t=0,d$ must have the same characteristic
exponents, i.e., $\gamma=\epsilon$.  Moreover, the exponents of $t=1$ must
be twice those of $t=0,d$, i.e., $1-\delta=2(1-\gamma)=2(1-\epsilon)$.
Either $R$ or~$S$ will be the degree-$4$ polynomial $4[t(2-t)-\frac12]^2$,
which maps $t=0,1,d$ to~$z=1$ {\rm(}$t=1$~with double multiplicity\/{\rm)}.
\end{enumerate}

The subcases of the `non-real' case~{\rm 2} are the following.

% change the numbering of the top-level enumerate environment
% to yield (2a),(2b),... instead of (1),(2),...
\def\labelenumi{{\em (\theenumi)}} \def\theenumi{2\alph{enumi}}

\begin{enumerate}
\item {\rm [}Equianharmonic {\rm(}equilateral triangle\,{\rm)} case.{\rm]}
$D=\frac12+ \ri\frac{\sqrt{3}}2$.  $q/\alpha\beta$~must equal the mean
of\/~$0,1,d$, and $t=0,1,d$ must have the same characteristic exponents,
i.e., $\gamma=\delta=\epsilon$.  Suppose $d=D$.  Then $q/\alpha\beta$ must
equal $\frac12+\ri\frac{\sqrt3}6$.  In~general, either $R$ or~$S$ will be
the degree-$3$ polynomial $\left[1-t/(\frac12+\ri\frac{\sqrt3}6)\right]^3$,
which maps $t=0,1,d$ to~$z=1$ and $t=q/\alpha\beta$ to~$z=0$ {\rm(}with
triple multiplicity\/{\rm)}.  There are special circumstances in which
$R$~may be sextic, which are listed separately, as subcase~{\rm 2d}.

\item $D=\frac12+\ri\frac{5\sqrt2}4$.  Suppose $d=D$.  Then $q/\alpha\beta$
must equal $\frac12+\ri\frac{\sqrt2}4$, $t=d$ must have characteristic
exponents $0,1/3$, i.e., $\epsilon=2/3$, and $t=0,1$ must have exponents
$0,1/2$, i.e., $\gamma=\delta=1/2$.  Either $R$ or~$S$ will be the
degree-$4$ polynomial
$\left[1-t/(\frac12+\ri\frac{5\sqrt2}4)\right]\left[1-t/(\frac12+\ri\frac{\sqrt2}4)\right]^3$,
which maps $t=d,q/\alpha\beta$ to~$z=0$ {\rm(}the latter with triple
multiplicity\/{\rm)} and $t=0,1$ to~$z=1$.

\item $D=\frac12+\ri\frac{11\sqrt{15}}{90}$.  Suppose $d=D$.  Then
$q/\alpha\beta$ must equal $\frac12+\ri\frac{\sqrt{15}}{18}$, $t=d$ must
have characteristic exponents $0,1/2$, i.e., $\epsilon=1/2$, and $t=0,1$
must have exponents $0,1/3$, i.e., $\gamma=\delta=2/3$.  Either $R$ or~$S$
will be the degree-$5$ polynomial
$At(t-1)\left[t-(\frac12+\ri\frac{\sqrt{15}}{18})\right]^3$, which maps
$t=0,1,q/\alpha\beta$ to~$z=0$ {\rm(}the last with triple
multiplicity\/{\rm)}.  The factor $A$ is chosen so that it maps $t=d$
to~$z=1$, as~well; explicitly, $A=-\ri\frac{2025\sqrt{15}}{64}$.

\item {\rm [}Special equianharmonic case.{\rm]}
$D=\frac12+\ri\frac{\sqrt{3}}2$.  $q/\alpha\beta$~must equal the mean
of\/~$0,1,d$, and $t=0,1,d$ must have characteristic exponents $0,1/3$,
i.e., $\gamma=\delta=\epsilon=2/3$.  Suppose $d=D$.  Then $q/\alpha\beta$
must equal $\frac12+\ri\frac{\sqrt3}6$.  Either $R$ or~$S$ will be the
degree-$6$ polynomial
$4\left\{[1-t/(\frac12+\ri\frac{\sqrt3}6)]^3-\frac12\right\}^2$, which maps
$t=0,1,d,q/\alpha\beta$ to~$z=1$ {\rm(}the last with triple
multiplicity\/{\rm)}.
\end{enumerate}
\end{thm}

% change the numbering of the top-level enumerate environment
% back to the default, not necessarily upright, (1),(2),...
\def\labelenumi{(\theenumi)}        \def\theenumi{\arabic{enumi}}

\begin{remarkaftertheorem}
The origin of the special harmonic and equianharmonic subcases is easy to
understand.  In subcase~1c, $t\mapsto R(t)$ or~$S(t)$ is the composition of
the quadratic map of subcase~1a with the map $z\mapsto4(z-\frac12)^2$.  In
subcase~2d, $t\mapsto R(t)$ or~$S(t)$ is similarly the composition of the
cubic map of subcase~2a with $z\mapsto4(z-\frac12)^2$.  In~both 1c~and~2d,
the further restrictions on exponents make possible the additional
quadratic transformation of~$z$, which transforms the hypergeometric
equation into itself~(see~\cite[\S\,3.1]{Andrews99} and~\cite{Erdelyi53}).
\end{remarkaftertheorem}

\begin{remarkaftertheorem}
\label{rem:rem2}
$R$ is determined uniquely by the choices enumerated in the theorem.  There
is a choice of subcase, a choice of~$d$ from the cross-ratio orbit of~$D$,
and a binary choice between $R$ and~$S$.  The final two choices amount to
choosing affine maps $A_1\in{\mathcal A}(\mathfrak{H})$
and~$A_2\in{\mathcal A}(\mathfrak{h})$, i.e., $A_2(z)=z$ or~$1-z$, which
precede and follow a canonical substitution.

In~the harmonic case~1a, in which the ${\mathcal A}(\mathfrak{H})$-orbit
includes three values of~$d$, there are accordingly $3\times2=6$
possibilities for~$R$; namely,
\begin{equation}
R=t^2,\,1-t^2;\quad (2t-1)^2,\,1-(2t-1)^2;\quad t(2-t),\,1-t(2-t),
\end{equation}
corresponding to $d=-1,-1;\frac12,\frac12;2,2$, respectively.  These are the
quadratic transformations of Kuiken~\cite{Kuiken79}.  In~the equianharmonic
case~2a, in which the orbit includes only two values of~$d$, there are
$2\times2=4$ possibilities; namely,
\begin{equation}
R=[1-t/(\tfrac12\pm \ri\tfrac{\sqrt3}6)]^3,
\qquad
1-[1-t/(\tfrac12\pm \ri\tfrac{\sqrt3}6)]^3,
\end{equation}
corresponding to $d=\frac12\pm \ri\frac{\sqrt{3}}2$.  The remaining
subcases, with the exception of 1c and~2d, correspond to generic
cross-ratio orbits: each value of~$D$ specifies six values of~$d$.  In~each
of those subcases, there are $6\times2=12$ possibilities.  So in~all, there
are 56~possibilities for~$R$.
\end{remarkaftertheorem}

\begin{remarkaftertheorem}
\label{rem:jumpahead}
The characteristic exponents of the singular points $z=0,1,\infty$
of~(\ref{eq:hyper}) can be computed from those of the singular points
$t=0,1,d,\infty$ of~(\ref{eq:Heun}), together with the formula for~$R$.
The computation relies on Proposition~\ref{thm:basicprop} below, which may
be summarized thus.  If $t=t_0$ is not a critical point of the map
$t\mapsto z=R(t)$, then the exponents of $z=R(t_0)$ will be the same as
those of~$t_0$.  If, on the other hand, $t=t_0$~is mapped to $z=z_0\defeq
R(t_0)$ with multiplicity~$k>1$, i.e., $t=t_0$ is a $k-1$-fold critical
point of~$R$ and $z=z_0$ is a critical value, then the exponents of $z_0$
will be $1/k$~times those of~$t_0$.

For example, in the harmonic case~1a, the map $t\mapsto z$ takes two of
$t=0,1,d$ to either $z=0$ or~$z=1$, and by examination, the coalesced point
is not a critical value of the map; so the characteristic exponents of
those two points are preserved, and must therefore be the same, as stated
in the theorem.  On~the other hand, the characteristic exponents of the
third point of the three, $t=d_0$, are necessarily halved when it is mapped
to $z=1$ or~$z=0$, since by examination, $R$~always has a simple critical
point at~$t=d_0$, i.e., $z\sim{\rm const} + C(t-d_0)^2$ for some
nonzero~$C$.  (These statements follow by considering the canonical $d=D$
case.)  So if $\delta_0$~denotes the parameter (out~of
$\gamma,\delta,\epsilon$) corresponding to~$t=d_0$, the characteristic
exponents of $z=1$ or~$z=0$ will be $0,(1-\delta_0)/2$.  $R$,~being a
quadratic polynomial, also has a simple critical point at~$t=\infty$, so
the characteristic exponents of~$z=\infty$ are one-half those
of~$t=\infty$, i.e., $\alpha/2,\beta/2$.  It~follows that in the harmonic
case, the Gauss parameters $(a,b;c)$ of the resulting hypergeometric
equation will be $(\alpha/2,\beta/2;(\delta_0+1)/2)$ or
$(\alpha/2,\beta/2;(\alpha+\beta-\delta_0+1)/2)$.

In~the equianharmonic case~2a, the map $t\mapsto z$ takes $t=0,1,d$ to
either $z=0$ or~$z=1$; and by examination, the coalesced point is not a
critical value of the map; so the characteristic exponents of those three
points are preserved, and must therefore be the same, as stated in the
theorem.  On~the other hand, at~$t=q/\alpha\beta$, which is mapped to $z=1$
or~$z=0$, $R$~has, by examination, a double critical point, i.e.,
$z\sim{\rm const}+ C(t-q/\alpha\beta)^3$ for some nonzero~$C$.  So~the
characteristic exponents of $z=1$ or~$z=0$, since $t=q/\alpha\beta$ is an
ordinary point of the Heun equation and effectively has characteristic
exponents~$0,1$, are~$0,1/3$.  $R$,~being a cubic polynomial, also has a
double critical point at~$t=\infty$, so the characteristic exponents
of~$z=\infty$ are one-third those of~$t=\infty$, i.e., $\alpha/3,\beta/3$.
It~follows that in the equianharmonic case, the parameters $(a,b;c)$ of the
resulting hypergeometric equation will be $(\alpha/3,\beta/3;2/3)$ or
$\left(\alpha/3,\beta/3;(\alpha+\beta+1)/3\right)$.
\end{remarkaftertheorem}

\begin{defn}
A rational map $R:\mathbb{CP}^1\to\mathbb{CP}^1$ is said to map the
characteristic exponents of the Heun equation~(\ref{eq:Heun}) to the
characteristic exponents of the hypergeometric
equation~($\mathfrak{h}$) if, for~all $t_0\in\mathbb{CP}^1$, the
exponents of $t=t_0$ according to the Heun equation, divided by the
multiplicity of $t_0\mapsto z_0\defeq R(t_0)$, equal the exponents
of~$z=z_0$ according to the hypergeometric equation.
\end{defn}

For example, if $t_0$ and~$z_0$ are both finite, this says that if $z\sim
z_0+C(t-t_0)^k$ to leading order, for some nonzero~$C$, then the exponents
of~$z=z_0$ must be those of~$t=t_0$, divided by~$k$.  If~$t=t_0$ is an
ordinary point of the Heun equation, then the exponents of~$z=z_0$ will
be~$0,1/k$.  This implies that if~$k>1$, $z=z_0$ must be one of the three
singular points of the hypergeometric equation.

\begin{prop}
\label{thm:basicprop}
A Heun equation of the form~{\rm(}\ref{eq:Heun}{\rm)} can be reduced to a
hypergeometric equation of the form~{\rm(}\ref{eq:hyper}{\rm)} by a
rational substitution $z=R(t)$ of its independent variable only~if $R$~maps
exponents to exponents.
\end{prop}

Proposition~\ref{thm:basicprop}, which was already used in
Remark~\ref{rem:jumpahead} above, is a special case of a basic fact in the
theory of the Riemann $P$-function: if~a rational change of the independent
variable transforms one Fuchsian equation to another, then the
characteristic exponents are transformed multiplicatively.  It~can be
proved by examining the effects of the change of variable on each local
(Frobenius) solution.  It~also follows immediately from
Lemma~\ref{thm:newlemma} below.

This lemma begins the study of {\em sufficient\/} conditions for the
existence of a Heun-to-hypergeometric transformation.  Finding them
requires care, since an accessory parameter is involved.  Performing the
substitution $z=R(t)$ explicitly is useful.  Substituting $z=R(t)$
into~(\ref{eq:hyper}) `pulls it back' (cf.~\cite{Kuiken79}) to
\begin{equation}
\label{eq:substituted}
\frac{\d^2 y}{\d t^2} + \left\{ -\frac{\ddot R}{\dot R} + \frac{\dot
R}{R(1-R)} [c-(a+b+1)R]\right\} \frac{\d y}{\d t} - \frac{ab\dot
R^2}{R(1-R)}\,y = 0.
\end{equation}
To save space here, $\d R/\!\d t,\d^2 R/\!\d t^2$ are written as $\dot
R,\ddot R$.

\begin{lem}
\label{thm:newlemma}
The coefficient of the $\d y/\!\d t$ term in the pulled-back hypergeometric
equation~{\rm(\ref{eq:substituted})}, which may be denoted~$W(t)$, equals
the coefficient of the $\d u/\!\d t$ term in the Heun
equation~{\rm(}\ref{eq:Heun}{\rm)}, i.e.,
$\gamma/t+\delta/(t-1)+\epsilon/(t-d)$, if~and only~if $R$~maps exponents
to exponents.  That~is, the transformation at~least partly `works' if~and
only~if $R$~maps exponents to exponents.
\end{lem}

\begin{pf}
This follows by elementary, if tedious calculations.  Suppose that $R$ maps
$t=t_0$ to $z=z_0\defeq R(t_0)$ with multiplicity~$k$, i.e., to~leading
order $R(t)\sim z_0+C(t-t_0)^k$; if $t_0$ and~$z_0$ are finite, that~is.
By~direct computation, the leading behavior of~$W$ at~$t=t_0$ is the
following.  In~the case when $t_0$~is finite, $W(t)\sim (1-k)(t-t_0)^{-1}$
if $z_0\neq0,1,\infty$; $[1-k(1-c)](t-t_0)^{-1}$ if $z_0=0$;
$[1-k(c-a-b)](t-t_0)^{-1}$ if $z_0=1$, and $[1-k(a+b)](t-t_0)^{-1}$ if
$z_0=\infty$.  In~the case when $t_0=\infty$, $W(t)\sim (1+k)t^{-1}$ if
$z_0\neq0,1,\infty$; $[1+k(1-c)]t^{-1}$ if $z_0=0$; $[1+k(c-a-b)]t^{-1}$ if
$z_0=1$, and $[1+k(a+b)]t^{-1}$ if $z_0=\infty$.

This may be restated as~follows.  At $t=t_0$, for finite~$t_0$, the leading
behavior of~$W$ is $W(t)\sim(1-k\eta)(t-t_0)^{-1}$, where $k$~is the
multiplicity of $t_0\mapsto z_0\defeq R(t_0)$ and $\eta$~is the sum of the
two characteristic exponents of the hypergeometric equation at~$z=z_0$;
if~the coefficient $1-k\eta$ equals zero then $W$~has no~pole at~$t=t_0$.
Likewise, the leading behavior of~$W$ at~$t=\infty$ is
$W(t)\sim(1+k\eta)t^{-1}$, where $k$~is the multiplicity of $\infty\mapsto
z_0\defeq R(\infty)$ and $\eta$ is the sum of the two exponents at~$z=z_0$;
if~the coefficient $1+k\eta$ equals zero then $W$~has a higher-order zero
at~$t=\infty$.

By the definition of `mapping exponents to exponents', it follows that the
leading behavior of~$W$ at $t=t_0$, for all $t_0$~finite, is of the form
$W(t)\sim(1-\eta')(t-t_0)^{-1}$, and also at $t_0=\infty$, is of the form
$W(t)\sim(1+\eta')t^{-1}$, where in both cases $\eta'$~is the sum of the
exponents of the Heun equation at~$t=t_0$, iff $R$~maps exponents to
exponents.

That is, the rational function~$W$ has leading behavior $\gamma t^{-1}$
at~$t=0$, $\delta(t-1)^{-1}$ at~$t=1$, $\epsilon(t-d)^{-1}$ at~$t=d$,
$(1+\alpha+\beta)t^{-1}=(\gamma+\delta+\epsilon)t^{-1}$ at~$t=\infty$, and
is regular at all~$t$ other than $0,1,d,\infty$, iff $R$~maps exponents to
exponents. \qed
\end{pf}

The following two propositions characterize when the pulled-back
hypergeometric equation~(\ref{eq:substituted}) is, in~fact, the Heun
equation~(\ref{eq:Heun}).  The first deals with Heun equations that are
trivial in the sense of Definition~\ref{def:triviality}, and will be used
in~\S\,\ref{sec:trivial}.  The second will be applied to prove
Theorem~\ref{thm:main}.

\begin{prop}
\label{thm:trivialprop}
A Heun equation of the form~{\rm(}\ref{eq:Heun}{\rm)}, which is trivial
{\rm(}i.e., $\alpha\beta=0$ and~$q=0${\rm)}, will be reduced to a
hypergeometric equation of the form~{\rm(}\ref{eq:hyper}{\rm)} by a
specified rational substitution $z=R(t)$ of its independent variable if~and
only~if $R$~maps exponents to exponents.
\end{prop}

\begin{pf}
The `if' half is new, and requires proof.  By Lemma~\ref{thm:newlemma}, the
coefficients of $\d y/\!\d t$ and $\d u/\!\d t$ agree iff $R$~maps
exponents to exponents, so it suffices to determine whether the
coefficients of $y$ and~$u$ agree.  But by triviality, the coefficient
of~$u$ in~{\rm (\ref{eq:Heun})} is zero.  Also, $t=\infty$ has zero as one
of its exponents, so all points $t\in\mathbb{CP}^1$ have zero as an
exponent.  By~the mapping of exponents to exponents, $z=\infty$ must also
have zero as an exponent, i.e., $ab=0$.  So~the coefficient of~$y$
in~(\ref{eq:substituted}) is also zero. \qed
\end{pf}

% change the numbering of the top-level enumerate environment
% to an upright version of (1),(2),...
\def\labelenumi{{\rm (\theenumi)}}        \def\theenumi{\arabic{enumi}}

\begin{prop}
\label{thm:mainproposition}
A Heun equation of the form~{\rm{\rm(}\ref{eq:Heun}{\rm)}}, which has four
singular points and is nontrivial {\rm(}i.e., $\alpha\beta\neq0$
or~$q\neq0${\rm)}, will be reduced to a hypergeometric equation of the
form~{\rm(}\ref{eq:hyper}{\rm)} by a specified rational substitution
$z=R(t)$ of its independent variable if~and only~if $R$~maps exponents to
exponents, and moreover, $R$~is a polynomial, $\alpha\beta\neq0$, and one
of the following two conditions on the normalized accessory parameter
$p\equiv q/\alpha\beta$ is satisfied.
\begin{enumerate}
\item $p$~equals one of\/ $0,1,d$.  Call this point~$d_0$, and the other
two singular points $d_1$~and~$d_2$.  In this case, $d_0$~must be a double
zero of $R$ or~$S$, and each of $d_1,d_2$ must be a simple zero of $R$
or~$S$.
\item $p$~does not equal any of\/ $0,1,d$.  In this case, each of\/ $0,1,d$
must be a simple zero of $R$ or~$S$, and $p$~must be a triple zero of
either $R$ or~$S$.
\end{enumerate}

% change the numbering of the top-level enumerate environment
% back to the default, not necessarily upright, (1),(2),...
\def\labelenumi{(\theenumi)}        \def\theenumi{\arabic{enumi}}

In both cases, $R$~and~$S$ must have no~additional simple zeroes or zeroes
of order greater than two.  Also, if\/ $1-c$ {\rm(}the nonzero exponent at
$z=0${\rm)} does not equal~$1/2$, then $R$~must have no~additional double
zeroes; and if $c-a-b$ {\rm(}the nonzero exponent at~$z=1${\rm)} does not
equal~$1/2$, then $S$~must have no~additional double zeroes.

Moreover, in both cases no~additional double zero, if~any, must be mapped
by~$R$ to the point\/ {\rm(}out~of $z=0,1${\rm)} to which $p$~is mapped.
{\rm(}So~additional double zeroes, if~any, must all be zeroes of~$R$, or
all be zeroes of~$S$.{\rm)}
\end{prop}

\begin{pf}
Like Proposition~\ref{thm:trivialprop}, this follows by comparing the
pulled-back hypergeometric equation~(\ref{eq:substituted}) to the Heun
equation~(\ref{eq:Heun}).  By Lemma~\ref{thm:newlemma}, the coefficients of
$\d y/\!\d t$ and $\d u/\!\d t$ agree iff $R$~maps exponents to exponents,
so it suffices to characterize when the coefficients of $y$ and~$u$ agree.

The coefficient of~$y$ in~(\ref{eq:substituted}) is to equal the
coefficient of~$u$ in~(\ref{eq:Heun}).  It~follows that $ab=0$ is
possible iff $\alpha\beta=0$ and~$q=0$, which is ruled~out by
nontriviality.  So $ab\neq0$, and equality of the coefficients can hold iff
\begin{equation}
\label{eq:cond}
U\equiv\frac{\d R/\!\d t}{R}\,\frac{\d S/\!\d t}{S} = 
\frac{-(\alpha\beta t - q)/ab}{t(t-1)(t-d)}
\equiv\frac{C_0}{t} + \frac{C_1}{t-1} + \frac{C_d}{t-d},
\end{equation}
where $S\equiv1-R$, and at~least two of $C_0,C_1,C_d\in\mathbb{C}$ are
nonzero.

Both $R^{-1}\d R/\!\d t$ and~$S^{-1}\d S/\!\d t$ are sums of terms of the form
$n(t-\lambda)^{-1}$, where $n$~is a nonzero integer and $\lambda$~is a zero
or a pole of $R$ or~$S$.  Poles are impossible, since $\lambda$~is a pole
of~$R$ iff $\lambda$~is a pole of~$S$, and there are no double poles on the
right-hand side of~(\ref{eq:cond}).  So $R$ must be a polynomial.

By examining the definition of~$U$ in~terms of $R$ and~$S$, one sees the
following is true of any $\lambda\in\mathbb{C}$\null: if $R$ or~$S$ has a
simple zero at~$t=\lambda$, then $U$~will have a simple pole
at~$t=\lambda$; if $R$ or~$S$ has a double zero at~$t=\lambda$, then
$U$~will have an ordinary point (non-zero, non-pole) at~$t=\lambda$, and if
$R$ or~$S$ has a zero of order $k>2$ at~$t=\lambda$, then $U$~will have a
zero of order $k-2$ at~$t=\lambda$.

Most of what follows is devoted to proving the `only~if' half of the
proposition in the light of these facts, by examining the consequences of
the equality~(\ref{eq:cond}).  In the final paragraph, the `if'~half will
be proved.

There are exactly three ways in which the equality~(\ref{eq:cond}) can
hold.

% By using `setcounter' appropriately, we start numbering with 0 here.
\begin{enumerate}
\setcounter{enumi}{-1}
\item $\alpha\beta=0$, but due to nontriviality, $q\neq0$.  $U$~has three
simple poles on~$\mathbb{C}$, at $t=0,1,d$.  It~has no other poles, and no
zeroes.  So each of~$0,1,d$ must be a simple zero of either $R$ or~$S$;
also, $R$~and~$S$ can have no~other simple zeroes, and no zeroes of
order~$k>2$.  Except for possible double zeroes, the zeroes of $R$ and~$S$
are determined.  The degree of~$R$ must equal the number of zeroes of~$R$,
and also equal the number of zeroes of~$S$, counting multiplicity.  But
irrespective of how many double zeroes are assigned to $R$ or~$S$, either
$R$ or~$S$ will have an odd number of zeroes, and the other an even number,
counting multiplicity.  So case~0 cannot occur.
\item $\alpha\beta\neq0$ and $\alpha\beta t-q$ is a nonzero multiple of
$t-d_0$, where $d_0=0$,~$1$, or~$d$, so exactly one of $C_0,C_1,C_d$ is
zero.  $U$~has two simple poles on~$\mathbb{C}$, at $t=d_1,d_2$ (the~two
singular points other than~$d_0$); it~has no other poles, and no zeroes.
So each of~$d_1,d_2$ must be a simple zero of either $R$ or~$S$; also,
$R$~and~$S$ can have no~other simple zeroes, and no zeroes of order~$k>2$.
Since by assumption ($\mathfrak{H}$)~has four singular points, each of
$0,1,d$ must be a singular point, so the coefficient of $\d y/\!\d t$
in~(\ref{eq:substituted}) must have a pole at~$t=d_0$, which implies that
$R$ or~$S$ must have a zero at~$d_0$ of the only remaining type: a~double
zero.
\item $\alpha\beta\neq0$ but $\alpha\beta t-q$ is not a multiple of
$t$,~$t-1$, or~$t-d$, so none of $C_0,C_1,C_d$ is zero.  $U$~has three
poles on~$\mathbb{C}$, and exactly one zero, at $t=p\equiv q/\alpha\beta$,
which is simple.  So each of~$0,1,d$ must be a simple zero of either $R$
or~$S$, and $q/\alpha\beta$~must be a triple zero of either $R$ or~$S$.
Also, $R$~and~$S$ can have no~other simple zeroes, and no other zeroes of
order~$k>2$.
\end{enumerate}

In cases 1,2, what remain to be determined are the (additional) double
zeroes of $R$ and~$S$, if~any.  That~is, it must be determined if any {\em
ordinary\/} point of the Heun equation can be mapped to $z=0$ or~$z=1$ with
double multiplicity.  But by Proposition~\ref{thm:basicprop}, $R$~can map
an ordinary point $t=t_0$ to~$z=0$ (resp.\ $z=1$) in this way only~if the
exponents of~$z=0$ (resp.\ $z=1$) are~$0,1/2$.

Suppose this occurs.  In case~1, if the exponents of~$t=p=d_0$ are
denoted~$0,\gamma_0$, the exponents of~$R(p)$ will be~$0,\gamma_0/2$, since
$t=p$ will be mapped with double multiplicity to~$z=R(p)$.  So if
$R(t_0)=R(p)$ then $\gamma_0$~must equal~$1$, which, since $q=\alpha\beta
d_0$, is ruled~out by the assumption that each of $0,1,d$, including~$d_0$,
is a singular point.  It~follows that in case~1, $R(t_0)\neq R(p)$.
A~related argument applies in case~2.  In case~2, the point~$p$ is an
ordinary point of~($\mathfrak{H}$), and a double critical point of the
$t\mapsto z$ map.  So as a singular point of~($\mathfrak{h}$), $R(p)$ must
have exponents~$0,1/3$.  It~follows that $R(t_0)=R(p)$ is impossible.

The `only if' half of the proposition has now been proved; the `if' half
remains.  Just as (\ref{eq:cond})~implies the stated conditions on~$R$, so
the stated conditions must be shown to imply~(\ref{eq:cond}).  But the
conditions on~$R$ are equivalent to the left and right-hand sides having
the same poles and zeroes, i.e., to their being the same up~to a constant
factor.  To~show the constant is unity, it~is enough to consider the limit
$t\to\infty$.  If $\deg R=n$, then $R^{-1}\d R/\!\d t\sim n/t$ and
$S^{-1}\d S/\!\d t\sim -n/t$, so~$U$, i.e., the left-hand side, has
asymptotic behavior~$-n^2/t^2$.  This will be the same as that of the
right-hand side if $(\alpha\beta)/(ab)=n^2$.  But $a=\alpha/n$ and
$b=\beta/n$ follow from the assumption that $R$~maps exponents to
exponents. \qed
\end{pf}

\noindent
Finally, we can prove the main theorem, with the aid of the polynomial
manipulation facilities of the {\sc Macsyma} computer algebra system.

\begin{pf*}{Proof (of Theorem~\ref{thm:main})}  By
Proposition~\ref{thm:mainproposition}, the preimages of $z=0,1$ under~$R$
must include $t=0,1,d$, and in case~2, $t=p\equiv q/\alpha\beta$.  They may
also include $l$~(additional) double zeroes of $R$ or of~$S$, which will be
denoted $t=a_1,\dotsc,a_l$.  Cases~1,2 of the theorem correspond to
cases~1,2 of the proposition, and the subcases of the theorem correspond to
distinct choices of~$l$.

Necessarily $\deg R=|R^{-1}(0)|=|R^{-1}(1)|$, where the inverse images are
defined as multisets rather than sets, to take multiplicity into account.
This places tight constraints on~$l$, since each of $0,1,d$ (and~$p$, in
case~2) may be assigned to either $R^{-1}(0)$ or~$R^{-1}(1)$, but by the
proposition, all of $a_1,\dotsc,a_l$ must be assigned, twice, to one or the
other.  In case~1, one of $0,1,d$ (denoted~$d_0$ in the proposition) has
multiplicity~$2$, and the other two (denoted~$d_1,d_2$) have
multiplicity~$1$.  It~follows that $0\leq l\leq 2$, with $\deg R=l+2$.  In
case~2, each of $0,1,d$ has multiplicity~$1$, and $p$~has multiplicity~$3$.
It~follows that $0\leq l\leq 3$, with $\deg R=l+3$.  Subcases are
as~follows.

\begin{enumerate}

% change the numbering of the top-level enumerate environment
% to yield an upright (1a),(1b),... instead of (1),(2),...
\def\labelenumi{{\rm (\theenumi)}} \def\theenumi{1\alph{enumi}}

\item Case~1, $l=0$, $\deg R=2$.  Necessarily $R^{-1}(0),R^{-1}(1)$ are
$\{d_0,d_0\},\allowbreak\{d_1,d_2\}$, or vice versa.  Without loss of
generality (w.l.o.g.), assume the latter, and assume $d_0=1$.  Then
$R^{-1}(0)=\{0,d\}$ and $R^{-1}(1)=\{1,1\}$, i.e., $S^{-1}(0)=\{1,1\}$ and
$S^{-1}(1)=\{0,d\}$.  Since $t=1$ is a double zero of~$S$, $S(t)=C(t-1)^2$
for some~$C$.  But $S(0)=1$, which implies $C=1$, and $S(d)=1$, which
implies~$d=2$.  So $S(t)=(t-1)^2$ and $R(t)=t(2-t)$.  Since $t=0,d$ are
both mapped singly to the singular point~$z=0$ by~$R$, their exponents must
be those of~$z=0$, and hence must be identical.

\item Case~1, $l=1$, $\deg R=3$.  Necessarily $R^{-1}(0),R^{-1}(1)$ are
$\{d_0,d_0,d_1\},\allowbreak\{d_2,a_1,a_1\}$, or vice versa.  W.l.o.g.,
assume the former, and also assume $d_0,d_1,d_2$ equal $1,d,0$,
respectively.  Then $R^{-1}(0)=\{1,1,d\}$ and $R^{-1}(1)=\{0,a_1,a_1\}$.
It~follows that $R(t)=(t-1)^2(1-t/d)$, where $d$~is determined by the
condition that the critical point of~$R$ other than~$t=1$ (i.e., $t=a_1$)
be mapped to~$1$.  Solving $\d R/\!\d t=0$ yields $a_1=(2d+1)/3$, and
substitution into $R(a_1)-1=0$ yields $d=4$ or~$-1/2$.  But the latter is
ruled~out: it would imply $a_1=0$, which is impossible.  So $d=4$ and
$a_1=3$.  Since $t=1,d$ are mapped to the singular point~$z=0$, doubly and
singly respectively, the exponents of $t=1$ must be twice those of~$z=0$,
and the exponents of $t=d$ must be the same as those of~$z=0$.

\item Case~1, $l=2$, $\deg R=4$.  Necessarily $R^{-1}(0)$ and~$R^{-1}(1)$
are $\{d_0,d_0,d_1,d_2\}$ and $\{a_1,a_1,a_2,a_2\}$, or vice versa.
W.l.o.g., assume the latter, and assume ${d_0=1}$.  Then
$R^{-1}(0)=\{a_1,a_1,a_2,a_2\}$ and $R^{-1}(1)=\{0,1,1,d\}$, i.e.,
$S^{-1}(0)=\{0,1,1,d\}$ and $S^{-1}(1)=\{a_1,a_1,a_2,a_2\}$.  So $S(t)$
equals $t(t-1)^2(t-d)$, where $d$~is determined by the condition that
$S$~must have two critical points other than $t=1$, i.e., $t=a_1,a_2$,
which are mapped by~$S$ to the same critical value (in~fact, to~$z=1$).
Computation yields $\d S/\!\d t=A(t-1)\left[4t^2-(3d+2)t+d\right]$, so
$a_1,a_2$ must be the roots of $4t^2-(3d+2)t+d$.  If~the corresponding
critical values are $Aw_1,Aw_2$, then $w_1,w_2$ are the roots of the
polynomial in~$w$ obtained by eliminating~$t$ between $w-S(t)/A$ and
$4t^2-(3d+2)t+d$.  Its discriminant turns~out to be proportional to
$(d-2)^2(9d^2-4d+4)^3$, so the criterion for equal values is that $d=2$ or
$9d^2-4d+4=0$.  But the latter can be ruled~out, since by examination it
would result in $a_1,a_2$ being equal.  So~$d=2$; $a_1,a_2=1\pm\sqrt2/2$;
and $S(t)=At(t-1)^2(t-2)$ with $A=-4$, so that $S(a_i)=1$.  Hence
$R(t)=4\left[t(2-t)-\frac12\right]^2$.  Since $t=0,d$ are mapped simply
to~$z=1$ and $t=1$~is mapped doubly, the exponents of~$t=0,d$ must be the
same, and double those of~$t=1$.

% change the numbering of the top-level enumerate environment
% to yield an upright (2a),(2b),...
\def\labelenumi{{\rm (\theenumi)}} \def\theenumi{2\alph{enumi}}
\setcounter{enumi}{0}

\item Case~2, $l=0$, $\deg R=3$.  Necessarily $R^{-1}(0),R^{-1}(1)$ are
$\{p,p,p\},\allowbreak\{0,1,d\}$, or vice versa.  W.l.o.g., assume the
former.  Then $R(t)=A(t-p)^3$ for some~$A$.  Since $t=0,1,d$ are to be
mapped singly to~$1$, they must be the vertices of an equilateral triangle,
with mean~$p$, so $A=-1/p^3$ and $R=(1-t/p)^3$.  W.l.o.g., take
$d=\frac12+\ri\frac{\sqrt3}2$, so $p=\frac12+\ri\frac{\sqrt3}6$.  The
exponents at $t=0,1,d$ must be equal, since they all equal the exponents
at~$z=1$.

\item Case~2, $l=1$, $\deg R=4$.  Assume w.l.o.g.\ that
$R^{-1}(0)=\{p,p,p,d\}$ and $R^{-1}(1)=\{0,1,a_1,a_1\}$, i.e.,
$S^{-1}(0)=\{0,1,a_1,a_1\}$ and $S^{-1}(1)=\{p,p,p,d\}$.  It~follows that
$R(t)=(1-t/d)(1-t/p)^3$, but to determine $d$ and~$p$, it~is best to focus
on~$S$.  Necessarily $S(t)=At(t-1)(t-a_1)^2$, and $p$~can be a triple zero
of~$R$ iff it is a double critical point of~$S$ as~well as~$R$.  The
condition that $S$~have a double critical point determines~$a_1$.  $\d
S/\!\d t=A(t-a_1)\left[4t^2-(3+2a_1)t+a_1\right]$, so the polynomial
$4t^2-(3+2a_1)t+a_1$ must have a double root.  Its discriminant is
$4a_1^2-4a_1+9$, which will equal zero iff $a_1=\frac12\pm\ri\sqrt{2}$.
The corresponding value of the double root, i.e., the mandatory value
of~$p$, is $\frac12\pm \ri\frac{\sqrt{2}}4$.  The requirement that $S$~map
$p$ to~$1$ implies~$A=1/p(p-1)(p-a_1)^2$.  $d$~is determined as the root
of~$R=1-S$ other than~$p$\/; some computation yields $\frac12\pm
\ri\frac{5\sqrt2}4$.  W.l.o.g.\ the `$\pm$' in the expressions for $p$
and~$d$ can be replaced by~`$+$'.  Since $t=p$ is an ordinary point and
$R$~maps $t=p$ triply to~$z=0$, $z=0$ must have exponents~$0,1/3$.  Since
$R$~maps $t=d$ simply to~$z=0$, $t=d$~must also have exponents~$0,1/3$.
Similarly, since $R$~maps the ordinary point $t=a_1$ doubly to~$z=1$,
$z=1$~must have exponents~$0,1/2$; so $t=0$ and~$t=1$, which are mapped
simply to~$z=1$, must also.

\item Case~2, $l=2$, $\deg R=5$.  Assume w.l.o.g.\ that
$R^{-1}(0)=\{p,p,p,0,1\}$ and $R^{-1}(1)=\{d,a_1,a_1,a_2,a_2\}$.  Then
$R(t)=At(t-1)(t-p)^3$, where $p$~is determined by $R$~having two critical
points other than $t=p$, i.e., $t=a_1,a_2$, which are mapped to the same
critical value (i.e., to ${z=1}$).  $\d R/\!\d
t=A(t-p)^2\left[5t^2-(2p+4)t+p\right]$, so $a_1,a_2$ must be the two roots
of $5t^2-(2p+4)t+p$.  If~the corresponding critical values are $Aw_1,Aw_2$,
then $w_1,w_2$ are the roots of the polynomial in~$w$ obtained by
eliminating~$t$ between $w-R(t)/A$ and $5t^2-(2p+4)t+p$.  Its discriminant
turns~out to be proportional to $(p^2-p+4)^3(27p^2-27p+8)^2$, so the
criterion for equal values is that $p^2-p+4=0$ or $27p^2-27p+8=0$.  But the
former can be ruled~out, since by examination it would result in $a_1,a_2$
being equal.  The latter is true iff $p=\frac12\pm
\ri\frac{\sqrt{15}}{18}$.  W.l.o.g., the plus sign may be used.  This
yields $a_1,a_2=\frac12 \pm \frac{2\sqrt3}9 + \ri\frac{\sqrt{15}}{90}$.
From the condition $R(a_i)=1$, it~follows that
$A=-\ri\frac{2025\sqrt{15}}{64}$.  $d$~is determined as the root of
$R(t)-1$ other than~$a_1,a_2$; computation yields
$d=\frac12+\frac{11\sqrt{15}}{90}$.  Since $t=p$ is an ordinary point
mapped triply to~$z=0$, $z=0$~must have exponents~$0,1/3$.  Similarly,
since $R$~maps the ordinary points~$t=a_i$ to~$z=1$, $z=1$~must have
exponents~$0,1/2$, so~$t=d$, which is mapped singly to~it, must also.

\item Case~2, $l=3$, $\deg R=6$.  Necessarily $R^{-1}(0)$ and~$R^{-1}(1)$
are $\{p,p,p,0,1,d\}$ and $\{a_1,a_1,a_2,a_2,a_3,a_3\}$, or vice versa.
W.l.o.g., assume the latter.  Then $R(t)=A(t-a_1)^2(t-a_2)^2(t-a_3)^2$ and
$S(t)=Bt(t-1)(t-d)(t-p)^3$.  Since $t=p$ is a triple zero of~$S$,
$R(t)\sim1-C(t-p)^3$ for some nonzero~$C$.  So~$\sqrt{R(t)}$, defined to
equal~$+1$ at~$t=p$, will have a similar Taylor series:
$\sqrt{A}(t-a_1)(t-a_2)(t-a_3)\sim1-C(t-p)^3/2$.  This is possible only~if
$a_1,a_2,a_3$ are the vertices of an equilateral triangle, and $p$~is their
mean.  It~follows that the roots of~$S$ other than~$t=p$, i.e., $t=0,1,d$,
are also the vertices of an equilateral triangle centered on~$p$.
W.l.o.g., choose $d=\frac12+\ri\frac{\sqrt3}2$
and~$p=\frac12+\ri\frac{\sqrt3}6$.  With a bit of algebra, $R$~can be
rewritten in the form given in the theorem.  Since $t=p$ is an ordinary
point and $R$~maps it triply to~$z=0$, $z=0$ must have exponents~$0,1/3$.
Since $R$~maps $t=0,1,d$ simply to~$z=0$, $t=0,1,d$~must also have
exponents~$0,1/3$.
% change the numbering of the top-level enumerate environment
% back to the default, not necessarily upright, (1),(2),...
\def\labelenumi{(\theenumi)}        \def\theenumi{\arabic{enumi}}
\end{enumerate}
The theorem is proved. \qed
\end{pf*}

Now that the main theorem is proved, we can proceed to derive explicit
Heun-to-hypergeometric reduction formul\ae.  Theorem~\ref{thm:culmination}
gives a necessary condition for the existence of a reduction, and
Theorem~\ref{thm:useful0} presents a representative list.

\begin{thm}
\label{thm:culmination}
Suppose the Heun equation\/~{\rm(\ref{eq:Heun})} has four singular points
and is nontrivial {\rm(}$\alpha\beta\neq0$ or~$q\neq0${\rm)}.  Then its
local solution $\Hl(d,q;\alpha,\beta,\gamma,\delta;t)$ can be reduced to a
hypergeometric function ${}_2F_1(a,b;c;z)$ by a formula of the type
$\Hl(t)={}_2F_1(R(t))$, with $R$ a rational function, only if its
parameters $d,q;\alpha,\beta$ satisfy $q=\alpha\beta p$, with $(d,p)$ equal
to one of the following\/ $23$~pairs.  If a reduction of this type exists,
$R(t)$~will be a polynomial of the stated degree.
\begin{displaymath}
\begin{array}{lll}
{\rm(1a)}&  [\deg R=2\hbox{ \rm{or} }4]. & (-1,0),\,(\frac12,\frac12),\,(2,1). \\
{\rm(1b)}&  [\deg R=3]. &
(-3,0),\,(-\frac13,0),\,(\frac14,\frac14),\,(\frac34,\frac34),\,(\frac43,1),\,(4,1). \\
{\rm(2a)}&  [\deg R=3\hbox{ \rm{or} }6]. & 
	(\frac12\pm\ri\frac{\sqrt{3}}2, \frac12\pm\ri\frac{\sqrt{3}}6). \\
{\rm(2b)}&  [\deg R=4]. & (\frac12\pm\ri\frac{5\sqrt2}4,\frac12\pm\ri\frac{\sqrt2}4), \\
&  & (\frac4{27}\pm\ri\frac{10\sqrt2}{27},\frac7{27}\pm\ri\frac{4\sqrt2}{27}),\,(\frac{23}{27}\pm\ri\frac{10\sqrt2}{27},\frac{20}{27}\pm\ri\frac{4\sqrt2}{27}). \\
{\rm(2c)}&  [\deg R=5]. &
(\frac12\pm\ri\frac{11\sqrt{15}}{90},\frac12\pm\ri\frac{\sqrt{15}}{18}), \\
&  & (\frac{135}{128}\pm\ri\frac{33\sqrt{15}}{128},\frac{95}{128}\pm\ri\frac{9\sqrt{15}}{128}),\,(-\frac{7}{128}\pm\ri\frac{33\sqrt{15}}{128},\frac{33}{128}\pm\ri\frac{9\sqrt{15}}{128}).
\end{array}
\end{displaymath}
\end{thm}
\begin{pf}
The five subcases are taken from Theorem~\ref{thm:main}, with (1c),(2d)
respectively subsumed into (1a),(2a).  Theorem~\ref{thm:main} supplies only
a single pair $(D,p_0)$ for each subcase.  By the remarks preceding that
theorem, $d$~may be any value on the cross-ratio orbit of~$D$.  There are
up~to $6$~distinct possibilities, which were given in~(\ref{eq:firstlist}).
The value of~$p$ associated to~$d$ is the preimage of~$p_0$ under the
corresponding affine map $t\mapsto A_1(t)$, given in~(\ref{eq:secondlist}).
\qed
\end{pf}

\begin{thm}
\label{thm:useful0}
Suppose a Heun equation has four singular points and is nontrivial
{\rm(}$\alpha\beta\neq0$ or~$q\neq0${\rm)}.  Then the only reductions of
its local Heun function~$\Hl$ to~${}_2F_1$ that can be performed by a
rational transformation of the independent variable involve polynomial
transformations of degrees $2$,~$3$, $4$, $5$, and\/~$6$.  There are seven
distinct types, each of which can exist only if $d$~lies on an appropriate
cross-ratio orbit.  The following list includes a representative reduction
of each type.  The ones with real~$d$ {\rm(}and $\deg R=2,3,4${\rm)}
include
\begin{subequations}
\begin{align}
\label{eq:generalharmonic}
\qquad&\Hl\left(2,\,\alpha\beta;\,\alpha,\,\beta,\,\gamma,\,\alpha+\beta-2\gamma+1;\,t\right)\\
\qquad&\qquad={}_2F_1\left(\tfrac{\alpha}2,\,\tfrac{\beta}2;\,\gamma;\,t(2-t)\right),
\nonumber\\[\jot]
\qquad&\Hl\left(4,\,\alpha\beta;\,\alpha,\,\beta,\,\tfrac12,\,\tfrac{2(\alpha+\beta)}3;\,t\right)\\
\qquad&\qquad={}_2F_1\left(\tfrac{\alpha}3,\,\tfrac{\beta}3;\,\tfrac12;\,1-(t-1)^2(1-t/4)\right),
\nonumber\\[\jot]
\label{eq:specialharmonic}
\qquad&\Hl\left(2,\,\alpha\beta;\,\alpha,\,\beta,\,\tfrac{\alpha+\beta+2}4,\,\tfrac{\alpha+\beta}2;\,t\right)\\
\qquad&\qquad={}_2F_1\left(\tfrac{\alpha}4,\,\tfrac{\beta}4;\,\tfrac{\alpha+\beta+2}4;\,1-4\bigl[t(2-t)-\tfrac12\bigr]^2\right),
\nonumber
\end{align}
\end{subequations}
and the ones with non-real~$d$ {\rm(}and $\deg R=3,4,5,6${\rm)} include
\begin{subequations}
\begin{align}
\label{eq:generalequianharmonic}
\qquad&\Hl\left(\tfrac12\pm\ri\tfrac{\sqrt3}{2},\,\alpha\beta(\tfrac12\pm\ri\tfrac{\sqrt3}6);\,\alpha,\,\beta,\,\tfrac{\alpha+\beta+1}3,\,\tfrac{\alpha+\beta+1}3;\,t\right)\\
\qquad&\qquad={}_2F_1\left(\tfrac{\alpha}3,\,\tfrac{\beta}3;\,\tfrac{\alpha+\beta+1}3;\,1-\bigl[1-t/(\tfrac12\pm\ri\tfrac{\sqrt3}6)\bigr]^3\right),
\nonumber\\[\jot]
\qquad&\Hl\left(\tfrac12\pm\ri\tfrac{5\sqrt2}4,\,\alpha({\tfrac23}-\alpha)(\tfrac12\pm\ri\tfrac{\sqrt2}4);\,\alpha,\,{\tfrac23}-\alpha,\,\tfrac12,\,\tfrac12;\,t\right)\\
\qquad&\qquad={}_2F_1\Bigl(\tfrac{\alpha}4,\,{\tfrac16}-\tfrac{\alpha}4;\,\tfrac12;\,
\nonumber\\[-4pt]
\qquad&\qquad\qquad\qquad 1-\bigl[1-t/(\tfrac12\pm\ri\tfrac{5\sqrt2}4)\bigr]\bigl[1-t/(\tfrac12\pm\ri\tfrac{\sqrt2}4)\bigr]^3\Bigr),
\nonumber
\displaybreak[0]
\\[\jot]
\label{eq:neednorm}
\qquad&\Hl\left(\tfrac12\pm\ri\tfrac{11\sqrt{15}}{90},\,
\alpha({\tfrac56}-\alpha)(\tfrac12\pm\ri\tfrac{\sqrt{15}}{18});\,
\alpha,\,{\tfrac56}-\alpha,\,{\tfrac23},\,{\tfrac23},\,t\right)\\
\qquad&\qquad={}_2F_1\Bigl(\tfrac{\alpha}5,\,{\tfrac16}-\tfrac{\alpha}5;\,{\tfrac23};\,
\nonumber\\[-4pt]
\qquad&\qquad\qquad\qquad (\mp\ri\tfrac{2025\sqrt{15}}{64})\,t(t-1)\bigl[t-(\tfrac12\pm\ri\tfrac{\sqrt{15}}{18})\bigr]^3\Bigr),
\nonumber\\[\jot]
\label{eq:specialequianharmonic}
\qquad&\Hl\left(\tfrac12\pm\ri\tfrac{\sqrt3}2,\,\alpha(1-\alpha)(\tfrac12\pm\ri\tfrac{\sqrt3}6);\,\alpha,\,1-\alpha,\,{\tfrac23},\,{\tfrac23},\,t\right)\\
\qquad&\qquad={}_2F_1\left(\tfrac{\alpha}6,\,\tfrac16-\tfrac{\alpha}6;\,{\tfrac23};\,1-4\left\{\bigl[1-t/(\tfrac12\pm\ri\tfrac{\sqrt3}6)\bigr]^3-\tfrac12\right\}^2\right).
\nonumber
\end{align}
\end{subequations}
In the preceding reductions, $\alpha,\beta,\gamma$ are free parameters.
Each of these equalities holds in a neighborhood of $t=0$ whenever the two
sides are defined, e.g., whenever the fifth argument of~$\Hl$ and the third
argument of~${}_2F_1$ are not equal to a nonpositive integer.
\end{thm}

\begin{remarkaftertheorem}
The equalities of Theorem~\ref{thm:useful0} hold even if the Heun
equation has fewer than four singular points, or is trivial; but in either
of those cases, additional reductions are possible.  For the trivial case,
see~\S\,\ref{sec:trivial}.
\end{remarkaftertheorem}

\begin{remarkaftertheorem}
The special harmonic reduction~(\ref{eq:specialharmonic}) is composite:
it~can be obtained from the case $\gamma=(\alpha+\beta+2)/4$ of the
reduction~(\ref{eq:generalharmonic}) by applying Gauss's quadratic
hypergeometric transformation~\cite[\S\,3.1]{Andrews99}
\begin{equation}
\label{eq:apply}
\quad{}_2F_1\bigl(a,\,b;\,(a+b+1)/2;\,z\bigr) = {}_2F_1\bigl(a/2,\,b/2;\,(a+b+1)/2;\,1-4(z-\tfrac12)^2\bigr)
\end{equation}
to the right-hand side.  The special equianharmonic
reduction~(\ref{eq:specialequianharmonic}) can be obtained in the same way
from the case $\beta=1-\alpha$ of the
reduction~(\ref{eq:generalequianharmonic}).

One might think that by applying (\ref{eq:apply}) to the right-hand sides
of the remaining reductions in
(\ref{eq:generalharmonic})--(\ref{eq:specialequianharmonic}), additional
composite reduction formul\ae\ could be generated.  However, there are only
a few cases in which it can be applied; and when it~can, it~imposes
conditions on the parameters of~$\Hl$ which require that the corresponding
Heun equation have fewer than four singular points.
\end{remarkaftertheorem}

\begin{pf}
$\Hl$ and~${}_2F_1$ are the local solutions of their respective equations
which belong to the exponent zero at~$t=0$ (resp.\ $z=0$), and are regular
and normalized to unity there.  So the theorem follows readily from
Theorem~\ref{thm:main}:
(\ref{eq:generalharmonic})--(\ref{eq:specialharmonic}) come from subcases
1a--1c, and
(\ref{eq:generalequianharmonic})--(\ref{eq:specialequianharmonic}) from
subcases 2a--2d.  In~each subcase, the Gauss parameters $(a,b;c)$
of~${}_2F_1$ are computed by first calculating the exponents at
$z=0,1,\infty$, in the way explained in Remark~\ref{rem:jumpahead}.
In~some subcases, the polynomial map supplied in Theorem~\ref{thm:main}
must be chosen to be $S=1-R$~rather than~$R$, due~to the need to map $t=0$
to $z=0$ rather than to~$z=1$, so~that the transformation will reduce~$\Hl$
to~${}_2F_1$, and not to another local solution of the hypergeometric
equation. \qed
\end{pf}

The list of Heun-to-hypergeometric reductions given in
Theorem~\ref{thm:useful0} is representative rather than exhaustive.  For
each subcase of Theorem~\ref{thm:main}, there is one reduction for each
allowed value of~$d$.  Each reduction on the above list came from choosing
$d=D$, but any other~$d$ on the cross-ratio orbit of~$D$ may be chosen.
The orbit is defined by $\triangle01d$ being one of the triangles (at~most
six) similar to~$\triangle01D$, i.e., by $\triangle01D$ being obtained
from~$\triangle01d$ by an affine transformation~$A_1\in{\mathcal
A}(\mathfrak{H})$.  So for any subcase of Theorem~\ref{thm:main} and choice
of~$d$, the appropriate polynomial map will be $z=A_2(R_1(A_1(t)))$, where
$A_1$~is constrained to map $\triangle01d$ to~$\triangle01D$ and is listed
in~(\ref{eq:secondlist}), $R_1$~is the polynomial map given in the subcase,
and $A_2\in{\mathcal A}(\mathfrak{h})$, i.e., $A_2(z)=z$ or~$1-z$, is
chosen so that $t=0$ is mapped to~$z=0$ rather than to~$z=1$.

As an example, consider the harmonic subcase~1a of Theorem~\ref{thm:main},
in which $D=2$, the cross-ratio orbit of~$D$ is~$\{-1,\frac12,2\}$, and the
polynomial map is $R_1(t)=t(2-t)$.  Choosing $d=D$ yields the
reduction~(\ref{eq:generalharmonic}).  Choosing $d=1-D=-1$ yields an
alternative reduction of~$\Hl$ to~${}_2F_1$, namely
\begin{equation}
\label{eq:extrathing}
\quad\Hl\bigl(-1,\,0;\,\alpha,\,\beta,\,\gamma,\,(\alpha+\beta-\gamma+1)/2;\,t\bigr)
={}_2F_1\left(\alpha/2,\,\beta/2;\,(\gamma+1)/2;\,t^2\right),
\end{equation}
in which $A_1(t)=1-t$ according to~(\ref{eq:secondlist}), and~$A_2(z)=1-z$.

It is not difficult to check that in~all, exactly $28$
Heun-to-hypergeometric reductions can be derived from
Theorem~\ref{thm:main}.  They exhibit the $23$~values of the pair
$(d,q/\alpha\beta)$ listed in Theorem~\ref{thm:culmination}.  Of~the~$28$,
eleven were given in Theorem~\ref{thm:useful0}, and
(\ref{eq:extrathing})~is a twelfth.  With the exception of the two
reductions with $(d,q/\alpha\beta)=(-1,0)$, one of which
is~(\ref{eq:extrathing}), the~$28$ split into pairs, each pair being
related by the identity~(\ref{eq:newguy1}), which takes $d$ to~$1/d$.

\section{A Generalization}
In applied mathematics, it is seldom the case that the four singular points
of an equation of Heun type are located at $0,1,d,\infty$.  But the main
theorem, Theorem~\ref{thm:main}, may readily be generalized.  Consider the
situation when three of the four have zero as a characteristic exponent,
since this may always be arranged by applying an F-homotopy.  There are two
cases of interest: either the singular points include~$\infty$ and each of
the finite singular points has zero as a characteristic exponent; or the
location of the singular points is unrestricted.  The latter includes the
former.  They have respective $P$-symbols
\begin{sizeequation}
{\small}
\label{eq:generalizedP}
P\left\{
\begin{array}{ccccc}
d_1&d_2&d_3&\infty& \\
0&0&0&\alpha& ;s \\
1-\gamma&1-\delta&1-\epsilon&\beta& 
\end{array}
\right\},
\quad
P\left\{
\begin{array}{ccccc}
d_1&d_2&d_3&d_4& \\
0&0&0&\alpha& ;s \\
1-\gamma&1-\delta&1-\epsilon&\beta& 
\end{array}
\right\}.
\end{sizeequation}
In the nomenclature of Ref.~\cite{Ronveaux95}, they are canonical cases of
the natural general-form and general-form Heun equations.  They are
transformed to the Heun equation~($\mathfrak{H}$) by the affine and
M\"obius transformations
\begin{equation}
\label{eq:respmaps}
t=\frac{s-d_1}{d_2-d_1},
\qquad
t=\frac{(s-d_1)(d_2-d_4)}{(d_2-d_1)(s-d_4)},
\end{equation}
respectively.  Each of the $P$-symbols~(\ref{eq:generalizedP}) is
accompanied by an accessory parameter.  The equation specified
by~(\ref{eq:generalizedP}a) can be written as
\begin{equation}
\label{eq:genHeun}
\quad\frac{\d^2 u}{\d s^2}
+ \left( \frac\gamma{s-d_1} + \frac\delta{s-d_2} + \frac\epsilon{s-d_3}
  \right)\frac{\d u}{\d s} + \frac{\alpha\beta s - q'}{(s-d_1)(s-d_2)(s-d_3)}\,u = 0,
\end{equation}
where $q'$ is the accessory parameter~\cite{Ronveaux95}.  The equation
specified by~(\ref{eq:generalizedP}b) with $d_4\neq\infty$ can be written
as
\begin{multline}
\label{eq:gen2Heun} 
\frac{\d^2 u}{\d s^2}
+ \left( \frac\gamma{s-d_1} + \frac\delta{s-d_2} + \frac\epsilon{s-d_3}
+ \frac{1-\alpha-\beta}{s-d_4}
\right)\frac{\d u}{\d s} \\
\qquad\qquad{}+ \frac{\alpha\beta\left.\left[{\displaystyle\prod_{i=1}^3(d_4-d_i)}\right]\right/(s-d_4) - q''}
{(s-d_1)(s-d_2)(s-d_3)(s-d_4)}\,u = 0,
\end{multline}
where $q''$ is the accessory parameter~\cite{Ronveaux95}.

The impediment to the generalization of Theorem~\ref{thm:main} to these two
equations is the specification of the cases that should be excluded due~to
their being `trivial', or having fewer than four singular points.  The
excluded cases should really be specified not in~terms of the {\em
ad~hoc\/} parameters $q'$ and~$q''$, but rather in an invariant way,
in~terms of an accessory parameter defined so~as to be invariant under
affine or M\"obius transformations, respectively.
Sch\"afke~\cite{Schafke83} has defined new accessory parameters of
second-order Fuchsian equations on~$\mathbb{CP}^1$ that are invariant under
affine transformations, but no~extension to general M\"obius
transformations seems to have been developed.

In the absence of an invariantly defined accessory parameter, an {\em
ad~hoc\/} approach will be followed.  It~is clear that
(\ref{eq:genHeun})~is trivial, i.e., can be transformed to a trivial Heun
equation by an affine transformation, iff $\alpha\beta=0$,~$q'=0$.  Also,
it will have fewer than four singular points if $\gamma=0$,~$q'=0$; or
$\delta=0$, $q'=\alpha\beta$; or $\epsilon=0$, $q'=\alpha\beta d$.
Likewise, it~is fairly clear that (\ref{eq:gen2Heun})~will be trivial,
i.e., can be transformed to a trivial Heun equation by a M\"obius
transformation, iff $\alpha\beta=0$,~$q''=0$.  The conditions on the
parameters for there to be a full set of singular points are, however, more
complicated.

The first generalization of Theorem~\ref{thm:main} is
Corollary~\ref{thm:gen1}, which follows from Theorem~\ref{thm:main} by
applying the affine transformation~(\ref{eq:respmaps}a).  It~mentions a
polynomial transformation, which is the composition of the $s\mapsto t$
affine transformation with the $t\mapsto z$ polynomial map of
Theorem~\ref{thm:main}.  To~avoid repetition, Corollary~\ref{thm:gen1}
simply cites Theorem~\ref{thm:main} for the necessary and sufficient
conditions on the exponent parameters and the accessory parameter.

\begin{cor}
\label{thm:gen1}
A natural general-form Heun equation of the canonical type\/
{\rm(\ref{eq:genHeun})}, which has four singular points and is nontrivial
{\rm(}i.e., $\alpha\beta\neq0$ or~$q'\neq0${\rm)}, can be reduced to a
hypergeometric equation of the form~{\rm{\rm(}\ref{eq:hyper}{\rm)}} by a
rational substitution $z=R(s)$ iff $\alpha\beta\neq0$, $R$~is a polynomial,
and the Heun equation satisfies the following conditions.

{\rm(i)}~$\triangle d_1d_2d_3$ must be similar to $\triangle 01D$, with
$D=2$ or $\frac12+\ri\tfrac{\sqrt3}2$, or $D=4$,
$\frac12+\ri\tfrac{5\sqrt2}4$, or~$\frac12+\ri\frac{11\sqrt{15}}{90}$.
That is, it must either be a degenerate triangle consisting of three
equally spaced collinear points {\rm(}the harmonic case\,{\rm)}, or be an
equilateral triangle {\rm(}the equianharmonic case\,{\rm)}, or be similar
to one of three other specified triangles, of which one is degenerate and
two are isosceles.  {\rm(ii)}~The exponent parameters
$\gamma,\delta,\epsilon$ must satisfy conditions that follow from the
corresponding subcases of Theorem\/~{\rm\ref{thm:main}}.  {\rm(iii)}~The
parameter~$q'$ must take a value that can be computed uniquely from the
parameters $\gamma,\delta,\epsilon$ and the choice of subcase.
\end{cor}

\begin{exampleaftertheorem}
In the harmonic case, the two endpoints of the degenerate triangle of
singular points $\triangle d_1d_2d_3$ must have equal exponent parameters,
and $q'$~must equal $\alpha\beta$ times the intermediate point.  In~this
case, $R$~will typically be a quadratic polynomial.  There are two
possibilities: $R$~will map the two endpoints to $z=0$ and the intermediate
point to $z=1$, or vice versa.  If~the characteristic exponents of the
intermediate point are twice those of the endpoints, then $R$~may be
quartic instead: the composition of either possible quadratic polynomial
with a subsequent ${z\mapsto 4(z-\frac12)^2}$ or $z\mapsto 4z(1-z)$ map.
\end{exampleaftertheorem}

\begin{exampleaftertheorem}
In~the equianharmonic case, all three exponent parameters
$\gamma,\delta,\epsilon$ must be equal, and the accessory parameter~$q'$
must equal $\alpha\beta$ times the mean of $d_1,d_2,d_3$.  In~this case,
$R$~will typically be a cubic polynomial.  There are two possibilities:
$R$~will map $d_1,d_2,d_3$ to $z=0$ and their mean to $z=1$, or vice versa.
If~the exponent parameters $\gamma,\delta,\epsilon$ equal~$2/3$, then
$R$~may be sextic instead: the composition of either possible cubic
polynomial with a subsequent $z\mapsto 4(z-\frac12)^2$ or $z\mapsto
4z(1-z)$ map.
\end{exampleaftertheorem}

The further generalization of Theorem~\ref{thm:main} is
Corollary~\ref{thm:gen2}, which follows from Theorem~\ref{thm:main} by
applying the M\"obius transformation~(\ref{eq:respmaps}b).  It~mentions a
rational substitution, which is the composition of the $s\mapsto t$
M\"obius transformation with the $t\mapsto z$ polynomial map of
Theorem~\ref{thm:main}.

\begin{cor}
\label{thm:gen2}
A general-form Heun equation of the canonical
type\/~{\rm(\ref{eq:gen2Heun})}, which has four singular points and is
nontrivial {\rm(}i.e., $\alpha\beta\neq0$ or~$q''\neq0${\rm)}, can be
reduced to a hypergeometric equation of the form~{\rm(}\ref{eq:hyper}{\rm)}
by a rational substitution $z=R(s)$ iff\/ $\alpha\beta\neq0$, and the Heun
equation satisfies the following conditions.

{\rm(i)}~The cross-ratio orbit of $\{d_1,d_2,d_3,d_4\}$ must be that of\/
$\{0,1,D,\infty\}$, where $D$~is one of the five values enumerated above.
That~is, it must be the harmonic orbit, the equianharmonic orbit, or one of
three specified generic orbits, one real and two non-real.  {\rm(ii)}~The
exponent parameters $\gamma,\delta,\epsilon$ must satisfy conditions that
follow from the corresponding subcases of Theorem\/~{\rm\ref{thm:main}}.
{\rm(iii)}~The parameter~$q''$ must take a value that can be computed
uniquely from the parameters $\gamma,\delta,\epsilon$ and the choice of
subcase.
\end{cor}

\begin{exampleaftertheorem}
Suppose $d_1,d_2,d_3,d_4$ form a harmonic quadruple, i.e., can be mapped by
a M\"obius transformation to the vertices of a square in~$\mathbb{C}$\null.
Moreover, suppose two of $d_1,d_2,d_3$ have the same characteristic
exponents, and are mapped to diagonally opposite vertices of the square.
That~is, of the three parameters $\gamma,\delta,\epsilon$, the two
corresponding to a diagonally opposite pair must be equal.  Then provided
$q''$~takes a value that can be computed from the other parameters, a
substitution~$R$ will exist.  It~will typically be a degree-$2$ rational
function, the only critical points of which are the third singular point
(out~of $d_1,d_2,d_3$) and~$d_4$.  Either $R$~will map the two
distinguished singular points to $z=1$ and the third singular point
to~$z=0$, or vice versa; and $d_4$ to~$z=\infty$.  In~the special case when
the exponents of the third point are twice those of the two distinguished
points, it~is possible for $R$ to be a degree-$4$ rational function.
\end{exampleaftertheorem}

\begin{exampleaftertheorem}
\label{example:2}
Suppose $d_1,d_2,d_3,d_4$ form an equianharmonic quadruple, i.e., can be
mapped by a M\"obius transformation to the vertices of a regular
tetrahedron in~$\mathbb{CP}^1$.  Moreover, suppose $d_1,d_2,d_3$ have the
same characteristic exponents, i.e., $\gamma=\delta=\epsilon$.  Then
provided $q''$~takes a value uniquely determined by the other parameters, a
substitution~$R$ will exist.  Typically, $R$ will be a degree-$3$ rational
function, the only critical points of which are the mean of $d_1,d_2,d_3$
with respect to~$d_4$, and~$d_4$.  Either $R$~will map $d_1,d_2,d_3$ to
$z=1$ and the mean of $d_1,d_2,d_3$ with respect to~$d_4$ to~$z=0$, or vice
versa; and $d_4$ to~$z=\infty$.  In~the special case when the exponents of
each of $d_1,d_2,d_3$ equal $0,1/3$, it~is possible for $R$ to be a
degree-$6$ rational function.
\end{exampleaftertheorem}

\begin{remarkaftertheorem}
In Example~\ref{example:2}, the concept of the mean of three points
in~$\mathbb{CP}^1$ with respect to a distinct fourth point was used.
A~projectively invariant definition is the following.  If $T$~is a M\"obius
transformation that takes $d_4$ ($\neq\nobreak d_1,d_2,d_3$) to the point at
infinity, the mean of $d_1,d_2,d_3$ with respect to~$d_4$ is the point that
would be mapped to the mean of $Td_1,Td_2,Td_3$ by~$T$.
\end{remarkaftertheorem}

\section{The Clarkson--Olver Transformation}
\label{sec:CO}
The reduction discovered by Clarkson and Olver~\cite{Clarkson96}, which
stimulated these investigations, turns~out to be a special case of the
equianharmonic Heun-to-hypergeometric reduction of~\S\,\ref{sec:main}.
Their reduction was originally given in a rather complicated form, which we
shall simplify.

Recall that the Weierstrass function $\wp(u)\equiv\wp(u;g_2,g_3)$ with
invariants $g_2,g_3\in\mathbb{C}$, which cannot both equal zero, has a
double pole at~$u=0$ and satisfies
\begin{align}
\label{eq:Wode} 
{\wp'}^2 &= 4\wp^3 - g_2\wp - g_3\\
&= 4(\wp - e_1)(\wp - e_2)(\wp - e_3).\nonumber
\end{align}
Here $e_1,e_2,e_3$, the zeroes of the defining cubic polynomial, are the
finite critical values of~$\wp$, the sum of which is zero; they are
required to be distinct.  $\wp$~is doubly periodic on~$\mathbb{C}$, with
periods denoted $2\omega,2\omega'$.  So it can be viewed as a function on
the torus $\mathbb{T}\defeq\mathbb{C}/{\mathcal L}$, where ${\mathcal
L}=2\omega\mathbb{Z}\oplus2\omega'\mathbb{Z}$ is the period lattice.  It
turns~out that the half-lattice $\{0,
\omega,\omega',\omega+\omega'\}+{\mathcal L}$ comprises the critical points
of~$\wp$.  The map $\wp:\mathbb{T}\to\mathbb{CP}^1$ is a double branched
cover of the Riemann sphere, but $\mathbb{T}$~is uniquely coordinatized by
the pair $(\wp,\wp')$.

The modular discriminant $\Delta\defeq g_2^3-27g_3^2\neq0$ is familiar from
elliptic function theory.  If $g_2,g_3\in\mathbb{R}$ and ${\Delta>0}$ (the
so-called real rectangular case, which predominates in applications),
$\omega,\omega'$ can be chosen to be real and imaginary, respectively.  If
$\Delta<0$ (the less familiar real rhombic case), they can be chosen to be
complex conjugates, so that the third critical point
$\omega_2\defeq\omega+\omega'$ is real.

Clarkson and Olver considered the Weierstrass-form Lam\'e equation
\begin{equation}
\label{eq:Lame}
\frac{\d^2\psi}{\d u^2} - \left[\ell(\ell+1)\wp(u) + B\right]\psi = 0,
\end{equation}
which is a Fuchsian equation on~$\mathbb{T}$ with exactly one singular
point (at~$(\wp,\wp')=(\infty,\infty)$) and a single accessory
parameter,~$B$.  [We~have altered their exponent parameter~$-36\sigma$ to
$\ell(\ell+1)$, to agree with the literature, and have added the accessory
parameter.]  In~particular, they considered the case $g_2=0$, $g_3\neq0$,
$B=0$.  They mapped $u\in\mathbb{T}$ to $z\in\mathbb{CP}^1$ by the formal
substitution
\begin{equation}
\label{eq:COsubst}
u = \frac{\ri}{(16g_3)^{1/6}}
\int^{(1-z)^{1/3}}\frac{\d\tau}{\sqrt{1-\tau^3}},
\end{equation}
and showed that the Lam\'e equation is reduced to
\begin{equation}
\label{eq:hyperCO}
z(1-z)\frac{\d^2\psi}{\d z^2} + \left(\frac12 - \frac76 z\right)
\frac{\d\psi}{\d z}
+\frac{\ell(\ell+1)}{36}\,\psi = 0.
\end{equation}
This is a hypergeometric equation with
$(a,b;c)=\left(-\ell/6,(\ell+1)/6;1/2\right)$.  It~has exponents $0,1/2$
at~$z=0$; $0,1/3$ at~$z=1$; and $-\ell/6,(\ell+1)/6$ at~$z=\infty$.

In elliptic function theory the case $g_2=0$, $g_3\neq0$ is called
equianharmonic, since the corresponding critical values $e_1,e_2,e_3$ are
the vertices of an equilateral triangle in~$\mathbb{C}$\null.  If, for
example, $g_3\in\mathbb{R}$, then $\Delta<0$; and by convention,
$e_1,e_2,e_3$ correspond to $\omega,\omega_2,\omega'$,
respectively. $e_1$~and~$e_3$ are complex conjugates, and $e_2$~is real.
The triangle $\triangle0\omega_2\omega'$ is also
equilateral~\cite[\S\,18.13]{Abramowitz65}.

So, what Clarkson and Olver considered was the {\em equianharmonic\/}
Lam\'e equation, the natural domain of definition of which is a
torus~$\mathbb T$ (i.e., a complex elliptic curve) with special symmetries.
For the Lam\'e equation~(\ref{eq:Lame}) to be viewed as a Heun equation
on~$\mathbb{CP}^1$, it~must be transformed by~$s=\wp(u)$ to its algebraic
form~\cite{Hille76}.  The algebraic form is
\begin{equation}
\label{eq:algLame}
\quad\frac{\d^2 \psi}{\d s^2}
+ \left( \frac{1/2}{s-e_1} + \frac{1/2}{s-e_2} + \frac{1/2}{s-e_3}
  \right)\frac{\d\psi}{\d s} + \frac{[-\ell(\ell+1)/4]s-B/4}{(s-e_1)(s-e_2)(s-e_3)}\,\psi = 0.
\end{equation}
This is a special case of~(\ref{eq:genHeun}), the canonical version of the
natural general-form Heun equation, with distinct finite singular points
$d_1,d_2,d_3 = e_1,e_2,e_3$.  Also, $\alpha,\beta=-\ell/2,(\ell+1)/2$,
$\gamma=\delta=\epsilon=1/2$, and~$q'=B/4$.  It~has characteristic
exponents $0,1/2$ at $s=e_1,e_2,e_3$, and $-\ell/2,(\ell+1)/2$
at~$s=\infty$.

Applying Corollary~\ref{thm:gen1} to~(\ref{eq:algLame}) yields the
following.

\begin{thm}
The algebraic-form Lam\'e equation~{\rm(\ref{eq:algLame})}, in the
equianharmonic case $g_2=0$, $g_3\neq0$, can be reduced when~$\ell(\ell+1)\neq0$
to a hypergeometric equation of the form~{\rm(}\ref{eq:hyper}{\rm)} by a
rational transformation $z=R(s)$ iff the accessory parameter~$B$ equals
zero.  If~this is the case, $R$~will necessarily be a cubic polynomial;
\begin{equation}
z=4s^3/g_3,\qquad z=1-4s^3/g_3
\end{equation}
will both work, and they are the only possibilities.
\end{thm}

\begin{pf}
If $\ell(\ell+1)\neq0$, the Heun equation~(\ref{eq:algLame}) is nontrivial
in the sense of Definition~\ref{def:triviality}, with four singular points;
by~(\ref{eq:Wode}), the~$e_i$ are the cube roots of~$g_3/4$, and are the
vertices of an equilateral triangle.  Since $\gamma=\delta=\epsilon$, the
equianharmonic case of Corollary~\ref{thm:gen1} applies, and no~other.

The mean of $e_1,e_2,e_3$ is zero.  So the polynomial $4s^3/g_3$ is the
cubic polynomial that maps each singular point to~$1$, and their mean to
zero; $1-4s^3/g_3$ does the reverse.  These are the only possibilities for
the map~$s\mapsto z$, since the sextic polynomials mentioned in the
equianharmonic case of Corollary~\ref{thm:gen1} can be employed only~if
$\gamma,\delta,\epsilon$ equal $2/3$, which is not the case here.  \qed
\end{pf}

{\em Remark\/}.  Corollary~\ref{thm:gen1} (equianharmonic case) also
applies to the equianharmonic algebraic-form Lam\'e equation with
$\ell(\ell+1)=0$, $B\neq0$, and guarantees it cannot be transformed to the
hypergeometric equation by any rational substitution; since in the sense
used above, this too is a nontrivial Heun equation.

\begin{cor}
\label{thm:cor}
The Weierstrass-form Lam\'e equation~{\rm(\ref{eq:Lame})}, in the
equianharmonic case $g_2=0$, $g_3\neq0$, can be reduced when $\ell(\ell+1)\neq0$
to a hypergeometric equation of the form~{\rm{\rm(}\ref{eq:hyper}{\rm)}} by
a substitution of the form $z=R\left(\wp(u)\right)$, where $R$~is rational,
iff the accessory parameter~$B$ equals zero.  In this case,
\begin{equation}
\label{eq:Wsubsts}
z=4\wp(u)^3/g_3,\qquad z=1-4\wp(u)^3/g_3
\end{equation}
will both work, and they are the only such substitutions.
\end{cor}

Applying the substitution $z=1-4\wp(u)^3/g_3$ to the Lam\'e
equation~(\ref{eq:Lame}) reduces it to the hypergeometric
equation~(\ref{eq:hyperCO}), as is readily verified.  (The other
substitution $z=4\wp(u)^3/g_3$ yields a closely related hypergeometric
equation, with the singular points $z=0,1$ interchanged.)  The
Clarkson--Olver substitution formula~(\ref{eq:COsubst}) contains a
multivalued elliptic integral, but it may be inverted with the aid
of~(\ref{eq:Wode}) to yield $z=1-4\wp(u)^3/g_3$.  So their transformation
fits into the framework of Corollary~\ref{thm:cor}.

The most noteworthy feature of the Clarkson--Olver transformation is that
it can be performed irrespective of the choice of exponent
parameter~$\ell$.  Only the accessory parameter~$B$ is restricted.  As~they
remark, when $\ell=1$,~$1/2$, $1/4$, or~$1/10$, it~is a classical result of
Schwarz that all solutions of the hypergeometric
equation~(\ref{eq:hyperCO}) are necessarily
algebraic~\cite[\S\,10.3]{Hille76}.  This implies that if~$B=0$, the same
is true of all solutions of the algebraic Lam\'e
equation~(\ref{eq:algLame}); which had previously been proved by
Baldassarri~\cite{Baldassarri81}, using rather different techniques.  But
irrespective of the choice of~$\ell$, the solutions of the $B=0$ Lam\'e
equation reduce to solutions of the hypergeometric equation.  This is quite
unlike the other known classes of exact solutions of the Lam\'e equation,
which restrict~$\ell$ to take values in a discrete
set~\cite[\S\,2.8.4]{Morales99}.  But it is typical of hypergeometric
reductions of the Heun equation.  As the theorems of~$\S\,\ref{sec:main}$
make clear, in~general it is possible to alter characteristic exponents
continuously, without affecting the existence of a reduction to the
hypergeometric equation.

It should be mentioned that the harmonic as~well as the equianharmonic case
of Corollary~\ref{thm:gen1} can be applied to the algebraic-form Lam\'e
equation.  One of the resulting quadratic transformations was recently
rediscovered by Ivanov~\cite{Ivanov2001}, in a heavily disguised form.  The
case of quadratic rather than cubic changes of the independent variable
will be considered elsewhere.

\section{The Seemingly Trivial Case $\alpha\beta=0$, $q=0$}
\label{sec:trivial}
If the Heun equation~(\ref{eq:Heun}) is trivial in the sense of
Definition~\ref{def:triviality}, it may be solved by quadratures.  A~basis
of solutions is
\begin{equation}
\quad u_1(t) = 1,\quad\,\,
u_2(t) = \int^t\exp\left[-\int^v\left(\frac\gamma w+\frac\delta{w-1} +
\frac\epsilon{w-d}\right)\,\d w\right]\,\d v.
\end{equation}
In the trivial limit, the local Heun function
$\Hl(d,q;\alpha,\beta,\gamma,\delta;t)$ degenerates to the former, and the
solution belonging to the exponent $1-\gamma$ at~$t=0$, denoted
$\widetilde\Hl(d,q;\alpha,\beta,\gamma,\delta;t)$ here, to the latter.
In~applications, explicit solutions, if~any, are what matter.  It~is
nonetheless interesting to examine under what circumstances a trivial Heun
equation can be reduced to a hypergeometric equation.  This question was
first considered by Kuiken~\cite{Kuiken79}.

The canonical polynomial substitutions of~\S\,\ref{sec:main} give rise to
many {\em nonpolynomial\/} rational reductions of trivial Heun equations to
hypergeometric equations, by composing with certain M\"obius
transformations.  To~understand why, recall that Theorem~\ref{thm:main}
characterized, up~to affine automorphisms of the two equations, the
polynomial substitutions that can reduce a nontrivial Heun equation to a
hypergeometric equation.  If $t\mapsto R_1(t)$ denotes a canonical
polynomial substitution, the full set of polynomial substitutions derived
from~it comprises all $t\mapsto A_2\left(R_1(A_1(t))\right)$, where
$A_1\in{\mathcal A}(\mathfrak{H})$~is an affine automorphism of the Heun
equation, which maps $\{0,1,d\}$ onto $\{0,1,D\}$, and $A_2\in{\mathcal
A}(\mathfrak{h})$~is an affine automorphism of the hypergeometric equation,
which maps $\{0,1\}$ onto~$\{0,1\}$.  (The only two possibilities are
$A_2(z)=z$ and $A_2(z)=1-z$.)

In the context of {\em nontrivial\/} Heun equations, M\"obius automorphisms
that are not affine could not be employed; essentially because, as
discussed in~\S\,\ref{subsec:auto}, moving the point at infinity would
require a compensating F-homotopy.  But in the trivial case no~such issue
arises: by Proposition~\ref{thm:trivialprop}, the Heun equation is reduced
to a hypergeometric equation by a rational substitution of its independent
variable, $z=R(t)$, iff the substitution maps exponents to exponents.  And
M\"obius transformations that are not affine certainly preserve exponents.

\begin{thm}
\label{thm:trivial}
A Heun equation of the form~{\rm(}\ref{eq:Heun}{\rm)}, which has four
singular points and is trivial {\rm(}i.e., $\alpha\beta=0$ and~$q=0${\rm)},
can be reduced to a hypergeometric equation of the
form~{\rm(}\ref{eq:hyper}{\rm)} by any rational substitution of the form
$z=M_2\left(R_1(M_1(t))\right)$, where $z=R_1(t)$~is a polynomial that maps
$\{0,1,D\}$ to $\{0,1\}$, listed {\rm(}along with~$D${\rm)} in one of the
seven subcases of Theorem\/~{\rm\ref{thm:main}}, and where $M_1\in{\mathcal
M}(\mathfrak{H})$, and $M_1\in{\mathcal M}(\mathfrak{h})$.  That~is,
$M_1$~maps $\{0,1,d,\infty\}$ onto $\{0,1,D,\infty\}$, and $M_2$~maps
$\{0,1,\infty\}$ onto $\{0,1,\infty\}$.  The necessary conditions on
characteristic exponents stated in Theorem\/~{\rm\ref{thm:main}} must be
satisfied, the conditions on exponents at specified values of~$t$ being
taken to refer to the exponents at the preimages of these points
under~$M_1$.
\end{thm}

\begin{remarkaftertheorem}
As in the derivation of the $\Hl(t)={}_2F_1(R(t))$ reduction formul\ae\
listed in Theorem~\ref{thm:useful0}, the Gauss parameters $(a,b;c)$ of the
resulting hypergeometric equation can be computed by first calculating the
exponents at~$z=R(t)=0,1,\infty$, using the mapping of exponents to
exponents.
\end{remarkaftertheorem}

The following example shows how such nonpolynomial rational substitutions
are constructed.  In the harmonic subcase~1a of Theorem~\ref{thm:main},
$D=2$ and the polynomial tranformation is $t\mapsto z=R_1(t)=t(2-t)$; the
necessary condition on exponents is that $t=0,d$ have identical exponents.
Consider $d=-1$, which is on the cross-ratio orbit of~$D$.
$M_1(t)=(t-1)/t$ can be chosen; also, let $M_2(z)=1/z$.  Then the
composition
\begin{equation}
z=R(t)\equiv M_2\left(R_1(M_1(t))\right)= t^2/(t^2-1)
\end{equation}
maps $t=0$ to~$z=0$ and $t=\infty$ to~$z=1$ (both with double
multiplicity), and $t=1,d$ to~$z=\infty$.  This substitution may be applied
to any trivial Heun equation with $d=-1$ and identical exponents at
$t=1,d$, i.e., with~${\delta=\epsilon}$.

In~this example, ${M}_1,{M}_2$ were selected with foresight, to ensure that
$R$~maps $t=0$ to~$z=0$.  This makes it possible to regard the substitution
as a reduction of $\Hl$ to~${}_2F_1$, or of $\widetilde\Hl$
to~$\widetilde{{}_2F_1}$.  By calculation of exponents, the reduction is
\begin{align}
\label{eq:explicitred}
&\widetilde\Hl\left(-1,\,0;\,0,\,\beta,\,\gamma,(1+\beta+\gamma)/2;\,t\right)\\
&\qquad=(-1)^{(\gamma-1)/2}\,
\widetilde{{}_2F_1}\left(0,\,(1-\beta+\gamma)/2;\,(1+\gamma)/2;\,t^2/(t^2-1)\right).\nonumber
\end{align}
The normalization factor $(-1)^{(\gamma-1)/2}$ is present because by
convention $\widetilde\Hl(t)\sim t^{1-\gamma}$ and
$\widetilde{{}_2F_1}(z)\sim z^{1-c}$ in a neighborhood of~$t=0$ (resp.\
$z=0$), where the principal branches are meant.  The corresponding
reduction of $\Hl$ to~${}_2F_1$ is trivially valid (both sides are constant
functions of~$t$, and equal unity).

Working out the number of rational substitutions $z=R(t)$ that may be
applied to trivial Heun equations, where $R$~is of the form
$M_2\circ R_1\circ M_1$,
%$M_2\left(R_1(M_1(\cdot))\right)$, 
is a useful exercise.  There are seven subcases of Theorem~\ref{thm:main},
i.e., choices for the polynomial~$R_1$.  Each subcase allows $d$~to be
chosen from an orbit consisting of $m$~cross-ratio values: $m=3$~in the
harmonic subcases 1a and~1c, $m=2$ in the equianharmonic subcases 2a
and~2d, and $m=6$ in the others.  In~any subcase, the $4!$~choices
for~$M_1$ are divided equally among the $m$~values of~$d$, and there are
also $3!$~choices for~$M_2$.  So each subcase yields $(4!/m)3!$ rational
substitutions for each value of~$d$, but not all are distinct.

To count {\em distinct\/} rational substitutions for each value of~$d$,
note the following.  $R$~will map $t=0,1,d,\infty$ to~$z=0,1,\infty$.  Each
of the subcases of Theorem~\ref{thm:main} has a `signature', specifying the
cardinalities of the inverse images of the points $0,1,\infty$.  For
example, case~1a has signature $2;1;1$, which means that of those three
points, one has two preimages and the other two have one.  (Order here is
not significant.)  In~all, subcases 1a,1b,2b,2c have signature $2;1;1$, and
the others have signature $3;1;0$.  By~inspection, the number of distinct
mappings of $t=0,1,d,\infty$ to~$z=0,1,\infty$ consistent with the
signature $2;1;1$ is~$36$, and the number consistent with $3;1;0$ is~$18$.

Kuiken~\cite{Kuiken79} supplies a useful list of the $36$~rational
substitutions arising from the harmonic subcase~1a, but states incorrectly
that they are the only rational substitutions that may be applied to a
trivial Heun equation.  Actually, subcases 1a--1c and 2a--2d give rise to
$36,36,18;\allowbreak18,36,36,18$ rational substitutions, respectively.  By
dividing by~$m$, it~follows that for each subcase, the number of distinct
rational substitutions per value of~$d$ is $12,6,6;\allowbreak9,6,6,9$.
Of~these, exactly one-third map $t=0$ to~$z=0$, rather than to $z=1$
or~$z=\infty$, and consequently yield reductions of $\Hl$ to~${}_2F_1$, or
of $\widetilde\Hl$ to~$\widetilde{{}_2F_1}$.  So for each subcase, the
number of such reductions per value of~$d$ is $4,2,2;\allowbreak3,2,2,3$.

For example, the four reductions with $d=-1$ that arise from the harmonic
subcase~1a are
\begin{subequations}
\begin{align}
\label{eq:firstred}
&\widetilde\Hl\left(-1,\,0;\,0,\,\beta,\,\gamma,\,(1+\beta-\gamma)/2;\,t\right)\\
&\qquad=\widetilde{{}_2F_1}\left(0,\,\beta/2;\,(1+\gamma)/2;\,t^2\right)\nonumber\\
\label{eq:secondred}
&\widetilde\Hl\left(-1,\,0;\,0,\,\beta,\,\gamma,(1+\beta+\gamma)/2;\,t\right)\\
&\qquad=(-1)^{(\gamma-1)/2}\,\widetilde{{}_2F_1}\left(0,\,(1-\beta+\gamma)/2;\,(1+\gamma)/2;\,t^2/(t^2-1)\right)\nonumber
\displaybreak[0]
\\
\label{eq:newred1}
&\widetilde\Hl\left(-1,\,0;\,0,\,\beta,\,1-\beta,\,\delta;\,t\right)\\
&\qquad=4^{-\beta}\widetilde{{}_2F_1}\left(0,\,(1-2\beta+\delta)/2;\,1-\beta;\,4t/(t+1)^2\right)\nonumber\\
\label{eq:newred2}
&\widetilde\Hl\left(-1,\,0;\,0,\,\beta,\,1-\beta,\,\delta,\,t\right)\\
&\qquad=(-4)^{-\beta}\,\widetilde{{}_2F_1}\left(0,\,(1-\delta)/2;\,1-\beta;\,-4t/(t-1)^2\right)\nonumber
\end{align}
\end{subequations}
The reduction~(\ref{eq:firstred}), which is the only one of the four in
which the degree-$2$ rational function $R$~is a polynomial, is simply the
trivial (i.e., $\alpha=0$) case of the quadratic
reduction~(\ref{eq:extrathing}).  The reduction~(\ref{eq:secondred}) was
derived above as~(\ref{eq:explicitred}), but (\ref{eq:newred1})
and~(\ref{eq:newred2}) are new.  They are related by composition with
$z\mapsto z/(z-1)$, i.e., by the involution in~${\mathcal M}(\mathfrak{h})$
that interchanges $z=1$ and~$z=\infty$.

Remarkably, many rational reductions of trivial Heun equations to the
hypergeometric equation are {\em not\/} derived from the polynomial
reductions of Theorem~\ref{thm:main}.  The following curious degree-$4$
reduction is an example.  The rational function
\begin{equation}
z=Q(t)\defeq1-\left(\frac{t-1-\ri}{t-1+\ri}\right)^4
=\frac{8\ri\,t(t-1)(t-2)}{(t-1+\ri)^4}
\end{equation}
takes $t=0,1,d\equiv2,\infty$ to~$z=0$; and $t=1\pm\ri$ to~$z=1,\infty$
(both with quadruple multiplicity).  By Proposition~\ref{thm:trivialprop},
a trivial Heun equation with $d=2$ will be reduced by~$Q$ to a
hypergeometric equation iff $Q$~maps exponents to exponents.  This
constrains the singular points $t=0,1,d,\infty$ to have the same exponents;
which by Fuchs's relation~(\ref{eq:Pconstraint}) is possible only if each
has exponents~$0,1/2$; which must also be the exponents of~$z=0$.  Also,
since $t=1\pm\ri$ are ordinary points of the Heun equation, with
exponents~$0,1$, the exponents of the hypergeometric equation at
$z=1,\infty$ must be~$0,1/4$.  It~follows that on the level of solutions,
the reduction is
\begin{equation}
\label{eq:tiring}
\widetilde\Hl(2,\,0;\,0,\,\tfrac12,\,\tfrac12,\,\tfrac12;\,t)
=(\ri/4)^{1/2}\,\widetilde{{}_2F_1}\left( 0,\,{\tfrac14};\,\tfrac12;\,
\frac{8\ri\,t(t-1)(t-2)}{(t-1+\ri)^4}\right),
\end{equation}
where the normalization factor $(\ri/4)^{1/2}$ follows from the known
behavior of the functions $\widetilde\Hl(t)$ and $\widetilde{{}_2F_1}(z)$
as~$t\to0$ and $z\to0$.  

Each of the preceding rational reductions of a trivial $\widetilde\Hl$ to a
$\widetilde{{}_2F_1}$ can be converted to a rational reduction of a {\em
nontrivial\/} $\Hl$ to a~${}_2F_1$, by using the definitions
(\ref{eq:tilde1}),(\ref{eq:newguy2}) of
$\widetilde{\Hl},\widetilde{{}_2F_1}$.  For example,
(\ref{eq:tiring})~implies
\begin{align}
\label{eq:finalreduction}
&\Hl(2,\,\tfrac34;\,\tfrac12,\,1,\,\tfrac32,\,\tfrac12;\,t)\\
&\quad=
(1-t)^{1/2}(1-t/2)^{1/2}\left[1-t/(1-\ri)\right]^{-2}
{}_2F_1\left(
\tfrac12,\,\tfrac34;\,\tfrac32;\,
\frac{8\ri t(t-1)(t-2)}{(t-1+\ri)^4}
\right).\nonumber
\end{align}
The equality~(\ref{eq:finalreduction}) holds in a neighborhood of~$t=0$
(both sides are real when $t$~is real and sufficiently small).  This
reduction is not related to the previously derived harmonic
reduction~(\ref{eq:generalharmonic}), in which $d=2$ also.  The pair
$(d,q/\alpha\beta)$ here equals $(2,\frac32)$, which is not listed in
Theorem~\ref{thm:culmination}.

The formula~(\ref{eq:finalreduction}) is a reduction of a nontrivial~$\Hl$
to a~${}_2F_1$, but of a more general type than has been considered in this
paper.  The underlying reduction of the Heun equation~(\ref{eq:Heun}) to
the hypergeometric equation~(\ref{eq:hyper}) includes a linear change of
the {\em dependent\/} variable, resembling a complicated F-homotopy, in
addition to a rational change of the independent variable.

\ack 
The author gratefully acknowledges the hospitality of the Texas
Institute for Computational and Applied Mathematics (TICAM).

% biblio: need Kulkarni paper, Borwein identity??  
% Mention in discussion section, esp. in
% connection with Heun<->Heun identities

% irreducible params, `specialization of irr params', q-transformations
% Erdelyi's inverse q-transformations and quadratic Heun <-> Heun identities

% Special solns
% -- Erdelyi's expansions and Craster's 2-term solutions
% -- Kimura 1970 has 2-term reps, Smirnov preprint has other

% algebraic rather than rational R?  Borwein identity

% Heun: connections and monodromy (see Slavyanov-Lay sec. 1.4)
% -- Central Two-Point Connection Problem, Sch\"afke-Schmidt
% -- Monodromy ref: Baider-Churchill

%\bibliographystyle{elsart-num} 
%\bibliography{general}

\begin{thebibliography}{10}

\bibitem{Abramowitz65}
M.~Abramowitz, I.~A. Stegun (Eds.), Handbook of Mathematical Functions, Dover,
  New York, 1965.

\bibitem{Andrews99}
G.~E. Andrews, R.~Askey, R.~Roy, Special Functions, Vol.~71 of Encyclopedia of
  Mathematics and Its Applications, Cambridge University Press, Cambridge, UK,
  1999.

\bibitem{Arscott81}
F.~M. Arscott, The land beyond {B}essel: A survey of higher special functions,
  in: Ordinary and Partial Differential Equations: Proceedings of the Sixth
  Dundee Conference, no. 846 in Lecture Notes in Mathematics, Springer-Verlag,
  1980, pp. 26--45.

\bibitem{Babister67}
A.~W. Babister, Transcendental Functions Satisfying Nonhomogeneous Linear
  Differential Equations, Macmillan, New York, 1967.

\bibitem{Baldassarri81}
F.~Baldassarri, On algebraic solutions of {L}am{\'e}'s differential equation,
  J.~Differential Equations 41~(1) (1981) 44--58.

\bibitem{Clarkson96}
P.~A. Clarkson, P.~J. Olver, Symmetry and the {Chazy} equation, J.~Differential
  Equations 124~(1) (1996) 225--246.

\bibitem{Craster98}
R.~V. Craster, V.~H. Ho{\`a}ng, Applications of {F}uchsian differential
  equations to free boundary problems, Proc. Roy.~Soc. London Ser.~A 454~(1972)
  (1998) 1241--1252.

\bibitem{Debosscher98}
A.~Debosscher, A unification of one-dimensional {Fokker}--{Planck} equations
  beyond hypergeometrics: Factorization solution method and eigenvalue schemes,
  Phys. Rev.~E 57~(1) (1998) 252--275.

\bibitem{Dwork84}
B.~Dwork, On {K}ummer's twenty-four solutions of the hypergeometric
  differential equation, Trans. Amer. Math. Soc. 285~(2) (1984) 497--521.

\bibitem{Erdelyi53}
A.~Erd{\'e}lyi (Ed.), Higher Transcendental Functions, McGraw--Hill, New York,
  1953--55, also known as The Bateman Manuscript Project.

\bibitem{Exton93}
H.~Exton, Solutions of {H}eun's equation, Bull. Soc. Math. Belg. S{\'e}r.~B
  45~(1) (1993) 49--57.

\bibitem{Grove85}
L.~C. Grove, C.~T. Benson, Finite Reflection Groups, 2nd Edition,
  Springer-Verlag, New York/Berlin, 1985.

\bibitem{Guttmann93}
A.~J. Guttmann, T.~Prellberg, Staircase polygons, elliptic integrals, {H}eun
  functions, and lattice {G}reen functions, Phys. Rev.~E 47~(4) (1993)
  R2233--R2236.

\bibitem{Heun1889}
K.~Heun, Zur {T}heorie der {R}iemann'schen {F}unctionen zweiter {O}rdnung mit
  vier {V}erzweigungspunkten, Math. Ann. 33 (1889) 161--179.

\bibitem{Hille76}
E.~Hille, Ordinary Differential Equations in the Complex Domain, Wiley, New
  York, 1976.

\bibitem{Ivanov2001}
P.~Ivanov, On {L}am{\'e}'s equation of a particular kind, J.~Phys.~A 34~(39)
  (2001) 8145--8150, available as arXiv:math-ph/0008008.

\bibitem{Joyce94}
G.~S. Joyce, On the cubic lattice {Green} functions, Proc. Roy.~Soc. London
  Ser.~A 445~(1924) (1994) 463--477.

\bibitem{Kuiken79}
K.~Kuiken, Heun's equation and the hypergeometric equation, {SIAM} J.~Math.
  Anal. 10~(3) (1979) 655--657.

\bibitem{Maier04}
R.~S. Maier, Algebraic solutions of the {L}am{\'e} equation, revisited,
  J.~Differential Equations 198~(1) (2004) 16--34, available as arXiv:math.CA/0206285.

\bibitem{Morales99}
J.~J. {Morales Ruiz}, Differential {G}alois Theory and Non-Integrability of
  {H}amiltonian Systems, Birk{\"a}user, Boston/Basel, 1999.

\bibitem{Poole36}
E.~G.~C. Poole, Linear Differential Equations, Oxford University Press, Oxford,
  1936.

\bibitem{Ronveaux95}
A.~Ronveaux (Ed.), Heun's Differential Equations, Oxford University Press,
  Oxford, 1995.

\bibitem{Schafke83} 
F.~W. Sch{\"a}fke, Zur (konfluenten) Fuchsschen
Differentialgleichungen 2. Ordnung, Analysis 3~(1--4) (1983) 101--122.

\bibitem{Schafke80a}
R.~Sch{\"a}fke, D.~Schmidt, The connection problem for general linear ordinary
  differential equations at two regular singular points with applications to
  the theory of special functions, {SIAM} J.~Math. Anal. 11~(5) (1980)
  848--862.

\bibitem{Schmitz94}
F.~Schmitz, B.~Fleck, On the propagation of linear 3-{D} hydrodynamic waves in
  plane non-isothermal atmospheres, Astron. Astrophys. Suppl. Ser. 106~(1)
  (1994) 129--139.

\bibitem{Snow52}
C.~Snow, Hypergeometric and {Legendre} Functions with Applications to Integral
  Equations of Potential Theory, 2nd Edition, no.~19 in Applied Mathematics
  Series, National Bureau of Standards, Washington,~DC, 1952.

\end{thebibliography}

\end{document}